\documentstyle[11pt]{article}
\newlength{\mytopmargin}
\newlength{\myleftmargin}
\newcommand{\la}{\lambda}

\newcommand{\de}{\delta}
\newcommand{\om}{\omega}
\def\rlx{\relax\leavevmode}
\def\inbar{\vrule height1.5ex width.4pt depth0pt}
\def\zz{\rlx\hbox{\small \sf Z\kern-.4em Z}}
\def\rr{\rlx\hbox{\scriptsize \rm I\kern-.18em R}}
\def\nn{\rlx\hbox{\rm I\kern-.18em N}}
\def\qq{\rlx\hbox{\,$\inbar\kern-.3em{\rm Q}$}}

\newcommand{\inn}[2]{\left\langle #1, #2 \right\rangle}
\newcommand{\qbin}[2]{\left[ \!\!\begin{array}{c} #1 \\#2 \end{array}\!\!
\right]}
\newcommand{\pbin}[2]{\left( \!\!\begin{array}{c} #1 \\#2 \end{array}\!\!
\right)}
\newcommand{\bin}[2]{({ {\scriptscriptstyle #1} \atop 
{\scriptscriptstyle #2} })}
\newtheorem{lemma}{Lemma}[section]
\newtheorem{df}[lemma]{Definition}

\newtheorem{cor}[lemma]{Corollary}
\newtheorem{prop}[lemma]{Proposition}

\begin{document}
\noindent
\begin{center}{ \Large\bf Multivariable Al-Salam\ \&\ Carlitz 
polynomials \\
associated with the type $A$ $q$-Dunkl kernel}
\end{center}
\vspace{5mm}

\noindent
\begin{center} T.H.~Baker\footnote{email: tbaker@maths.mu.oz.au} and
P.J.~Forrester\footnote{email: matpjf@maths.mu.oz.au}\\[2mm]
{\it Department of Mathematics, University of Melbourne, \\
Parkville, Victoria 3052, Australia}
\end{center}
\vspace{.5cm}

\begin{quote} The Al-Salam\ \&\ Carlitz polynomials are $q$-generalizations
of the classical Hermite polynomials. Multivariable generalizations of
these polynomials are introduced via a generating function
involving a multivariable hypergeometric function which is the
$q$-analogue of the type-$A$ Dunkl integral kernel.
An eigenoperator is established for these polynomials and this is used
to prove orthogonality with respect to a certain Jackson integral
inner product. This inner product is normalized by deriving a
$q$-analogue of the Mehta integral, and the corresponding normalization
of the multivariable  Al-Salam\ \&\ Carlitz polynomials is derived
from a Pieri-type formula. Various other special properties of the
polynomials are also presented, including their relationship
to the shifted Macdonald polynomials and the big-$q$ Jacobi
polynomials.

\end{quote}

\section{Introduction}
The advent of Macdonald polynomials \cite{mac} has opened the
way to a theory of hypergeometric functions in which the
Macdonald polynomials are the underlying basis. Progress in this
direction has been made by Macdonald \cite{macunp1}, Kaneko
\cite{kaneko96a} and Lassalle \cite{lass97a}. 
Present in the work of
 Macdonald and  Lassalle is the hypergeometric function
\begin{equation} \label{cano}
{}_0{\cal F}_0(x;y;q,t) := \sum_{\kappa} {t^{b(\kappa)} \over
h_\kappa'(q,t) P(1,t,\dots,t^{n-1};q,t)} P_\kappa(x;q,t) P_\kappa(y;q,t).
\end{equation}
where $x$ (similarly $y$)
denotes the $n$-tuple $(x_1,\dots,x_n)$,
 the sum is
over all partitions of length $n$ and $P_\kappa$
 denotes the Macdonald polynomial normalized so
that the coefficient of the leading monomial is unity, while Kaneko
defines the hypergeometric function
\begin{equation} \label{kan}
{}_0{\psi}_0(x;y;q,t) := \sum_{\kappa} {(-1)^{|\kappa|}
q^{b(\kappa')} \over
h_\kappa'(q,t) 
P(1,t,\dots,t^{n-1};q,t)} P_\kappa(x;q,t) P_\kappa(y;q,t).
\end{equation}
The quantities $h_\kappa'(q,t)$ and $b(\kappa)$ in (\ref{cano}) and
(\ref{kan}) are defined by
\begin{equation}\label{db}
h_\kappa'(q,t) := \prod_{(i,j)\in \kappa} (1 - q^{\kappa_i - j +1}
t^{\kappa_j'-i}), \quad b(\kappa) = \sum_{j=1}^n (j-1) \kappa_j
\end{equation}
and $\kappa'$ in (\ref{kan}) denotes
 the partition conjugate to $\kappa$.

Using (\ref{db}), the explicit formulas \cite{mac}
\begin{equation}\label{ef}
t^{-b(\kappa)} 
 P(1,t,\dots,t^{n-1};q,t) = {\prod_{(i,j) \in \kappa}
(1 - q^{j-1}t^{n-i+1}) \over h_\kappa(q,t)}, \quad
 h_\kappa(q,t) = \prod_{(i,j) \in \kappa}(1 - q^{\kappa_i - j}
t^{\kappa_j'-i+1})
\end{equation}
and the facts that $P_\kappa(x;q,t) = P_\kappa(x;q^{-1},t^{-1})$,
$c^{|\kappa|}P_\kappa(x;q,t) = P_\kappa(cx;q,t)$, we see that
the hypergeometric functions (\ref{cano}) and (\ref{kan}) are
related by
\begin{equation}\label{re}
{}_0{\cal F}_0(x;y;q^{-1},t^{-1}) = 
{}_0{\psi}_0(x;t^{n-1}qy;q,t)
\end{equation} 
Furthermore, with $t=q^{1/\alpha}$, in the limit $q \to 1$ (\ref{cano}) with
$y$ replaced by $(1-q)y$ and (\ref{kan}) with 
$y$ replaced by $-(1-q)y$ become equal:
\begin{equation}\label{qt1}
\lim_{q \to 1} {}_0{\cal F}_0(x;(1-q)y;q,q^{1/\alpha})=
\lim_{q \to 1} {}_0{\psi}_0(x;-(1-q)y;q,t) = 
{}_0{\cal F}_0^{(\alpha)}(x;y)
\end{equation}
where
\begin{equation}\label{0f0}
{}_0^{}{\cal F}_0^{(\alpha)}(x;y) := \sum_{\kappa} \frac{\alpha^{|\kappa|}}
{d'_{\kappa}}
\frac{P^{(\alpha)}_{\kappa}(x)\,P^{(\alpha)}_{\kappa}(y)}
{P^{(\alpha)}_{\kappa}(1^n)}
\end{equation}
and the limit is to be taken term-wise in the series.
In (\ref{0f0}), $P^{(\alpha)}_{\kappa}(x)$ denotes the 
Jack polynomial normalized so that the coefficient of the leading
monomial is unity and
$$
d'_{\kappa} =  \prod_{(i,j)\in \kappa} \Big ( \alpha (\kappa_i - j
+1) + (\kappa_j' - i) \Big ).
$$

The hypergeometric function ${}_0{\cal F}_0(x;y)$ has a number
of noteworthy features. 
It was first observed by Lassalle (see \cite{forr96a} for a 
published version of most of the results of this work) that this
function is the explicit
realization of the symmetric type $A$ integral kernel occurring
in the work of Dunkl \cite{dunkl91a,dunkl92a}. As such it occurs in
multidimensional generalizations of the Laplace, Hankel and
Fourier transforms, as well as in the study of
the type $A$ Calogero-Sutherland quantum many body
system (see e.g.~\cite{forr96a}). Moreover, associated with
this kernel are a class of multivariable polynomials
generalizing the classical Hermite polynomials 
\cite{macunp1,lass96a,forr96a,vandiej96b,roesler97a}. It is the objective 
of this article to investigate the $q$-analogues of these polynomials
associated with
${}_0{\cal F}_0(x;y;q,t)$ and
${}_0{\psi}_0(x;y;q,t)$.

After reviewing some of the salient properties of the single-variable
Al-Salam\&Carlitz polynomials $U_n^{(a)}$, Section 2 contains the
definition of the multivariable Al-Salam\&Carlitz polynomials
$U_\kappa^{(a)}$ in terms of a generating function involving
${}_0{\cal F}_0(x;y;q,t)$. In Section 3 these polynomials are related
to the Macdonald polynomials by a $q$-exponential operator formula,
and they are also shown to satisy an eigenvalue equation. Three
different forms of the eigenoperator are given, one of which is
manifestly Hermitian with respect to a particular Jackson integral
inner product, thus establishing an orthogonality relation.

In Section 4 special properties of the 
 multivariable Al-Salam\&Carlitz polynomials are established. One of
these properties is a Pieri-type formula expressing
$e_1U_\kappa^{(a)}$, where $e_1$ denotes the first elementary
symmetric function, in terms of a linear combination of
$U_\lambda^{(a)}$. This formula is used to calculate the
normalization of $U_\kappa^{(a)}$ with respect to the Jackson
integral inner product, and the normalization in turn is used in
the derivation of some integral representations.
Section 5 contains a discussion on the relationship between
the multivariable Al-Salam\&Carlitz polynomials and the big
$q$-Jacobi polynomials introduced recently by Stokman
\cite{stok97a}.
We conclude
with two Appendices. The first is an expansion formula required
in calculating the normalization of the inner product, while in
the second a determinant formula for the $U_\kappa^{(a)}$ in terms of
the $U_n^{(a)}$ is given in the special case $t=q$.

\section{Generalization of the polynomials of Al-Salam and
Carlitz}
\setcounter{equation}{0}
\subsection{The case $n=1$}
For $n=1$ we have ${}_0{\cal F}_0^{(\alpha)}(x;y)=e^{xy}$. Replacing
$y$ by $2y$ and multiplying by $e^{-x^2}$ gives the generating
function for the classical Hermite polynomials:
\begin{equation}\label{gfch}
e^{-x^2} e^{2xy} = \sum_{n=0}^\infty {H_n(y) x^n \over n!}.
\end{equation}
We recall that the classical Hermite polynomials are eigenfunctions
of the differential operator
\begin{equation}
{d^2 \over d y^2} - 2 y {d \over dy}
\end{equation}
and have the orthogonality
\begin{equation}\label{cho}
\int_{-\infty}^\infty dy \,e^{-y^2} H_m(y) H_n(y) = 
\sqrt{\pi} 2^n n! \delta_{m,n}.
\end{equation}

In the $q$-case for $n=1$
\begin{equation}\label{e}
{}_0^{}{\cal F}_0^{(\alpha)}(x;y;q,t) =e_q(xy) := \sum_{n=0}^\infty {(xy)^n \over 
(q;q)_n} = {1 \over (xy;q)_\infty}
\end{equation}
\begin{equation}\label{E}
{}_0{\psi}_0(x;y;q,t) = E_q(-xy) := \sum_{n=0}^\infty
{(-1)^n q^{n(n-1)/2} (xy)^n \over (q;q)_n } = (xy;q)_\infty,
\end{equation}
($(x;q)_n$, $(x;q)_\infty$ have their usual meaning) where the
final equalities are valid for $|q| < 1$. The $q$-generalizations
of (\ref{gfch}) associated with (\ref{e}) and (\ref{E}) 
have been given by Al-Salam and Carlitz \cite{al65a}.

Consider first (\ref{e}). Al-Salam and Carlitz defined a family
of polynomials $\{U_n^{(a)}(y;q)\}$, $(a<0)$ by the generating
function
\begin{equation}\label{gfu}
\rho_a(x;q) e_q(xy) = \sum_{n=0}^\infty U_n^{(a)}(y;q)
{x^n \over (q;q)_n}, \quad \rho_a(x;q) := 
(x;q)_\infty (ax;q)_\infty = E_q(-x) E_q(-ax)
\end{equation}
and established the orthogonality
\begin{equation} \label{ou}
\int_a^1 w_U^{(a)}(x;q) U_m^{(a)}(x;q) U_n^{(a)}(x;q) \,d_qx
= (1-q) (-a)^n q^{n(n-1)/2} (q;q)_n \delta_{m,n},
\end{equation}
where
\begin{equation}\label{wu}
 w_U^{(a)}(x;q) := {(qx;q)_\infty ({qx \over a};q)_\infty \over
(q;q)_\infty (a;q)_\infty ({q \over a};q)_\infty}
\end{equation}
\begin{equation}\label{defqi}
\int_a^1 f(x) \, d_qx := (1-q) \left(
\sum_{n=0}^\infty f(q^n) q^n - a\sum_{n=0}^\infty f(aq^n) q^n \right),
\quad (a <0)
\end{equation}
In the case $a=-1$ the polynomials $U_n^{(a)}(x;q)$ are referred to
as the discrete $q$-Hermite polynomials \cite{gasperbk}, and they
reduce to the classical Hermite polynomials in the $q \to 1$
limit:
\begin{equation}\label{1dh}
\lim_{q \to 1} (1 - q)^{-n/2} U_n^{(-1)}(x(1-q)^{1/2};q) = 2^{-n/2}
H_n({x \over \sqrt{2}}).
\end{equation}
For general $a$, consideration of the three term recurrence for 
$U_n^{(a)}$ (eq.~(\ref{3term}) below) shows
\begin{equation}
\lim_{q \to 1} (1 - q)^{-n/2} U_n^{(-q^r)}(x(1-q)^{1/2};q) = 
2^{-n/2}
H_n({x-r \over \sqrt{2}}).
\end{equation}

In relation to (\ref{E}), Al-Salam and Carlitz introduced the 
polynomials $\{V_n^{(a)}(x;q)\}$ via the generating function
\begin{equation}\label{gfv}
{1 \over \rho_a(x;q)} E_q(-xy) = \sum_{n=0}^\infty V_n^{(a)}(y;q)
{(-1)^n q^{n(n-1)/2} x^n \over (q;q)_n}.
\end{equation}
Note from (\ref{gfu}) and (\ref{e}) that
\begin{equation}
{1 \over \rho_a(x;q)} = 
{1 \over (x;q)_\infty (ax;q)_\infty} = e_q(x) e_q(ax).
\end{equation}
These polynomials were shown to satisfy the orthogonality
\begin{equation}\label{ov}
\int_1^\infty w_V(x;q) V_m^{(a)}(x;q) V_n^{(a)}(x;q)\,d_qx =
(1-q) a^m (q;q)_m q^{-m^2} \delta_{m,n}
\end{equation}
where
\begin{equation}
w_V(x;q) = {(q ;q)_\infty ({1 \over a};q)_\infty 
(qa;q)_\infty \over {(x;q)'}_\infty
({x \over a};q)_\infty}
\end{equation}
(the dash in ${(x;q)'}_\infty$ denotes that any factor which 
vanishes for a particular choice of $x$ is to be deleted) with
\begin{equation}
\int_1^\infty f(x) d_qx := (1-q) \sum_{n=0}^\infty f(q^{-n})
q^{-n}
\end{equation}
Furthermore it was established in \cite{al65a} that the polynomials
$\{U_n^{(a)}(x;q)\}$ and $\{V_n^{(a)}(x;q)\}$ are simply related:
\begin{equation}\label{uv}
V_n^{(a)}(x;q) = U_n^{(a)}(x;q^{-1}).
\end{equation}
This can be seen by first noting from (\ref{e}), (\ref{E}) and
(\ref{re}) that
\begin{equation}\label{eE}
E_{q^{-1}}(-x) = e_q(qx)
\end{equation}
and using this in (\ref{gfv}) with $q$ replaced by $q^{-1}$ to
obtain the generating function
\begin{equation}
\rho_a(qx;q) e_q(qxy) = \sum_{n=0}^\infty V_n^{(a)} (y;q^{-1})
{(qx)^n \over (q;q)_n}
\end{equation}
Comparison with (\ref{gfu}) implies (\ref{uv}).

Using (\ref{uv}) we can check that the orthogonality (\ref{ov}) is
equivalent to (\ref{ou}) with $q$ replaced by $q^{-1}$. This is 
done by first generalizing $w_U^{(a)}$ to include an auxiliary
parameter $\mu$:
\begin{equation}
w_U^{(a)}(x;q,\mu) := {E_q(-\mu qx) E_q(-\mu {q x \over a}) \over
E_q(-\mu q) E_q(-a) E_q(-\mu {q \over a})}
\end{equation}\label{gwu}
With this definition (\ref{eE}) gives
\begin{equation}\label{wui}
w_U^{(a)}(x;q^{-1},\mu) = {e_q(\mu x) e_q(\mu {x \over a}) \over
e_q(\mu) e_q(qa) e_q({\mu \over a})}.
\end{equation}
Now, comparing (\ref{gwu}) and (\ref{wu}), we see from (\ref{defqi})
that 
\begin{equation}\label{lqi}
{1 \over 1 - q} 
\int_a^1 w_U^{(a)}(x;q) f(x) \, d_qx :=  
\sum_{n=0}^\infty \lim_{\mu \to 1}w_U^{(a)}(q^n;q,\mu)
f(q^n) q^n - a\sum_{n=0}^\infty 
\lim_{\mu \to 1}w_U^{(a)}(aq^n;q,\mu)f(aq^n) q^n.
\end{equation}
Next we set $a=-q^p$ for some fixed $p \in\zz^+$, and replace 
$q$ by $q^{-1}$ in (\ref{lqi}). From (\ref{wui}) and the product
formula in (\ref{e}), for $x=-aq^n$ ($n=0,1,\dots$),
$
\lim_{\mu \to 1} w_U^{(a)}(x;q^{-1},\mu) = 0
$
while for $x=q^n$, $n=0,1,\dots$,
$$
\lim_{\mu \to 1} w_U^{(a)}(x;q^{-1},\mu) =
w_V^{(a)} (x;q).
$$
Substituting in (\ref{lqi}) gives
\begin{equation}\label{qiuv}
{1 \over 1 - q}
\int_a^1 w_U^{(a)}(x;q) f(x) \, d_qx \Big |_{q \mapsto q^{-1}} =
\sum_{n=0}^\infty w_V^{(a)}(q^{-n};q^{-1}) f(q^{-n}) q^{-n},
\end{equation}
which  is the required inter-relationship.
A simple lemma of Stembridge \cite{stem88a} show that this formula,
proved valid for $a=-q^p$, must remain true for all $a$ such
that both sides are defined.

\subsection{Generating function for general $n$}
For general $n$ the generating function for the classical
Hermite polynomials
(\ref{gfch}) has the generalization \cite{lass96a,forr96a}
\begin{equation}\label{gfherm}
e^{-x_1^2}\cdots e^{-x_n^2} {}_0^{}{\cal F}_0^{(\alpha)}(x;2y) = 
\sum_{\kappa} {\alpha^{|\kappa|} H_\kappa(y;\alpha) 
P_\kappa^{(\alpha)}(x)
\over d_\kappa'},
\end{equation}
where the polynomials $\{H_\kappa(y;\alpha)\}$ are referred to
as the generalized Hermite polynomials. These
polynomials have an expansion
in terms of $\{P_\sigma^{(\alpha)}(y)\}$ of the form
\begin{equation}\label{hp}
H_\kappa(y;\alpha) = {2^{|\kappa|} \over P_\kappa^{(\alpha)}(1^n)} 
P_\kappa^{(\alpha)}(y) + \sum_{|\mu| < |\kappa|} c_{\kappa \mu}
P_\mu^{(\alpha)}(y).
\end{equation}
With
\begin{equation}
D_0:= \sum_{i=1}^n \frac{\partial^2}{\partial x_i^2} +
\frac{2}{\alpha}\sum_{i\neq j} \frac{1}{x_i-x_j}\frac{\partial}
{\partial x_i}, \qquad E_1 :=  
\sum_{i=1}^n x_i\,\frac{\partial}{\partial x_i},
\end{equation}
the generalized Hermite polynomials
are also related to the Jack polynomials via the exponential
operator formula \cite[eq.(3.21)]{lass96a,forr96a}
\begin{equation}\label{expop}
{2^{|\kappa|} \over P_\kappa^{(\alpha)}(1^n)}
\exp \Big ( - {1 \over 4} D_0 \Big )\, P_\kappa^{(\alpha)}(x)
=  H_\kappa(x;\alpha) 
\end{equation}
and they are
eigenfunctions of the operator
\begin{equation}\label{heo}
\tilde{H}^{(H)} := D_0 - 2 E_1 \label{hop}.
\end{equation}
They also
satisfy the orthogonality
\begin{equation}
\int_{(-\infty, \infty)^n}
  H_\kappa(y;\alpha)  H_\sigma(y;\alpha) \,
d\mu^{(H)}(y) = {2^{|\kappa|}
d_\kappa' {\cal N}_0^{(H)} \over  P_\kappa^{(\alpha)}(1^n)}
\delta_{\kappa, \sigma}
\end{equation}
where 
\begin{equation}\label{mi}
d\mu^{(H)}(y) := \prod_{j=1}^n e^{-y_j^2} \prod_{1 \le j < k \le n}
|y_j - y_k|^{2/\alpha} \,dy_1 \dots dy_n,
\end{equation}
\begin{equation}
{\cal N}_0^{(H)} :=  \int_{(-\infty, \infty)^n} d\mu^{(H)}(y) =
2^{-n(n-1)/2\alpha}
\pi^{n/2} \prod_{j=0}^{n-1}{\Gamma(1+(j+1)/\alpha) \over
 \Gamma (1 + 1/\alpha )}.
\end{equation}

In light of this result, the single variable polynomials of Al-Salam
and Carlitz suggest an $n$-dimensional generalization. Thus
we introduce two sets of multivariable polynomials $\{U_\kappa^{(a)}
(y;q,t)\}$ and $\{V_\kappa^{(a)}(y;q,t)\}$ by the generating
functions
\begin{equation}\label{gfmu}
\rho_a(x_1;q) \cdots \rho_a(x_n;q) \,
{}_0{\cal F}_0(x;y;q,t)
= \sum_{\kappa} {t^{b(\kappa)}U_\kappa^{(a)}(y;q,t)
P_\kappa(x;q,t) \over h_\kappa'(q,t) P_\kappa(1,t,\dots,t^{n-1};q,t)}
\end{equation}
and
\begin{equation}\label{gfmv}
{1 \over \rho_a(t^{-(n-1)}x_1;q) \cdots \rho_a(t^{-(n-1)}x_n;q)} \,
{}_0{\psi}_0(x;y;q,t)
= \sum_{\kappa} {(-1)^{|\kappa|} q^{b(\kappa')}V_\kappa^{(a)}(y;q,t)
P_\kappa(x;q,t) \over h_\kappa'(q,t) P_\kappa(1,t,\dots,t^{n-1};q,t)}.
\end{equation}
With these definitions the polynomials 
$U_\kappa^{(a)}$ and $V_\kappa^{(a)}$ are simply related.
This is seen by replacing $q,t$ by $q^{-1},t^{-1}$ in (\ref{gfmu}). Use of
(\ref{eE}) and (\ref{re}) and comparison with (\ref{gfmv})
shows
\begin{equation}\label{inter}
V_\kappa^{(a)}(y;q,t) = U_\kappa^{(a)}(y;q^{-1},t^{-1})
\end{equation}
in analogy with (\ref{uv}). Also, it follows from the
generating function (\ref{gfmu}) that the polynomials 
$U_\kappa^{(a)}$ have an expansion in terms of the 
Macdonald polynomials of the form
\begin{equation}\label{up}
U_\kappa^{(a)}(x;q,t) = P_\kappa(x;q,t) +
\sum_{|\nu| < |\kappa|} a_{\kappa \nu}P_\nu(x;q,t) 
\end{equation}
in analogy with (\ref{hp}).

The generalized Hermite polynomials can be reclaimed from the
polynomials $U_\kappa^{(-1)}$ by a limiting procedure analogous
to (\ref{1dh}). To see this put $q=t^\alpha$ and $a=-1$ in (\ref{gfmu}),
and note from (\ref{gfu}) and (\ref{E}) that $\rho_1(x;q) = (x^2;q^2)_\infty$.
Then replace $x$ by $(1-q)^{1/2}x$ and $y$ by $(1-q)^{1/2}y$, and consider 
the limit $q \to 1$. Now ${}_0 {\cal F}_0((1-q)^{1/2}x;(1-q)^{1/2}y;q,t)$
$= {}_0 {\cal F}_0(x;(1-q)y;q,t)$ so we can take the limit in this term
by using (\ref{qt1}), while it follows from (\ref{E}) that
$(x^2(1-q);q^2)_\infty \to e^{-x^2/2}$ as $q \to 1$. On the r.h.s.~of
(\ref{gfmu}) we note that in this limit
$$
{P_\kappa((1-q)^{1/2}x;q,q^{1/\alpha}) 
\over P_\kappa(1,t,\dots,t^{n-1};q,q^{1/\alpha})}
= (1-q)^{|\kappa|/2} {P_\kappa(x;q,q^{1/\alpha}) \over 
 P_\kappa(1,t,\dots,t^{n-1};q,q^{1/\alpha})}  \sim (1-q)^{|\kappa|/2}
{P_\kappa^{(\alpha)}(x) \over P_\kappa^{(\alpha)}(1^n)},
$$
where $P_\kappa^{(\alpha)}(x)$ denotes the Jack polynomial, while
$(1-q)^{|\kappa|/2} h_\kappa'(q,q^{1/\alpha}) \to d_\kappa'/\alpha^{|\kappa|}$.
Thus we have
$$
e^{-x_1^2/2} \cdots e^{-x_n^2/2} {}_0^{}{\cal F}_0^{(\alpha)}(x;y) =
\sum_\kappa {{\alpha}^{|\kappa|} \lim_{q \to 1}
(1 - q)^{-|\kappa|/2} U^{(-1)}_\kappa((1-q)^{1/2}y;q,q^{1/\alpha})
P_\kappa^{(\alpha)}(x) \over d_\kappa' P_\kappa^{(\alpha)}(1^n)}.
$$
Comparison with (\ref{gfherm}) gives the multidimensional
generalization of (\ref{1dh}),
\begin{equation}
\lim_{q \to 1} (1-q)^{-|\kappa|/2} U^{(-1)}_\kappa
((1-q)^{1/2}x;q,q^{1/\alpha}) =
2^{-|\kappa|/2}  P_\kappa^{(\alpha)}(1^n) H_\kappa({x \over \sqrt{2}};
\alpha).
\end{equation}
\section{Exponential operator formula, 
eigenvalue equation and Hermiticity}
\setcounter{equation}{0}
To further develop the theory of the polynomials
$\{U_\kappa^{(a)}(x;q,t)\}$, we must revise the definitions of
operators $Y_i^{\pm 1}$ and 
$D_i$, which play a role in the theory of non-symmetric Macdonald polynomials.
These operators in turn are defined in terms of certain operators $T_i$, 
$\omega$ appearing in the theory of type $A$ affine Hecke algebras.
The Demazure-Lustig operators $T_i$ are defined by
\begin{equation} \label{pucallpa}
T_i := t + \frac{tx_i-x_{i+1}}{x_i - x_{i+1}}\left(s_i -1 \right)
\hspace{2cm} i=1,\ldots,n-1 
\end{equation}
while
\begin{equation}\label{omega}
\om := s_{n-1}\cdots s_2\,s_1\tau_1 =
s_{n-1}\cdots s_i\tau_i s_{i-1}\cdots s_1 .
\end{equation}
For future reference we note that the operators $T_i$ and
$\omega$ have the properties
\begin{eqnarray}
T_i^{-1}\,x_{i+1} =t^{-1}x_i\,T_i \hspace*{2cm} &
T_i^{-1}\,x_i = x_{i+1}\,T_i^{-1} + (t^{-1}-1)x_i  \nonumber\\
T_i\,x_i = t x_{i+1}\,T_i^{-1} \hspace*{2cm} &
T_i\,x_{i+1} = x_i\,T_i +(t-1)x_{i+1} \label{id.2}\\
\om\,x_i = qx_n\om \hspace{2cm} & \om\,x_{i+1} = x_i\,\om 
 \nonumber
\end{eqnarray}
valid for $1\leq i\leq n-1$.

Using the definitions (\ref{pucallpa}) and (\ref{omega}) the operators
$Y_i$ \cite{cher91a,mac95} are defined by 
\begin{equation} \label{yops}
Y_i = t^{-n+i}\;T_i\cdots T_{n-1}\;\om\;T_1^{-1}\cdots T_{i-1}^{-1},
\hspace{2cm} 1\leq i\leq n.
\end{equation}
They mutually commute and have the following relations with the  operators $T_i$,
\begin{equation}\label{yt1}
T_i\,Y_{i+1}\,T_i = t\,Y_i \hspace{2cm} [T_i, Y_j] =0, \quad j\neq i, i+1
\end{equation}
The $q$-analogue of the type $A$ Dunkl operators were defined in \cite{forr97b}
by means of
\begin{equation}\label{qdunkl-def}
D_i := x_i^{-1}\;\left( 1-t^{n-1}\left[ 1 +(t^{-1}-1)\sum_{j=i+1}^n
t^{j-i}\,T_{ij}^{-1} \right] Y_i \right)
\end{equation}
where for $i<j$,
\begin{eqnarray}
T_{ij}^{-1} &:=& T_i^{-1}\;T_{i+1}^{-1}\;\cdots\;T_{j-2}^{-1}\;T_{j-1}^{-1}
\;T_{j-2}^{-1}\; \cdots\; T_i^{-1} \nonumber   \\
&=& T_{j-1}^{-1}\; T_{j-2}^{-1}\;\cdots\;T_{i+1}^{-1}\;T_i^{-1}\;
T_{i+1}^{-1}\;\cdots\;T_{j-1}^{-1}
\end{eqnarray}

\subsection{$q$-exponential operator formula}
{}From Theorem 5.2(c) and Proposition 5.4 of ref.~\cite{forr97b}, it
follows that for any symmetric function $f$ analytic in the neighbourhood
of the origin,
\begin{equation}\label{fF}
f(D_1^{(x)},\dots,D_n^{(x)}) {}_0{\cal F}_0(x;y;q,t) = 
f(y_1,\dots,y_n) {}_0{\cal F}_0(x;y;q,t)
\end{equation}
where the $D_i$ are the operators (\ref{qdunkl-def}) acting on the
$x$-variables only. By choosing $f(x_1,\dots,x_n) =
\rho_a(x_1;q) \cdots \rho_a(x_n;q)$ and using the generating function
(\ref{gfmu}) we see immediately that
\begin{equation}\label{expopu}
\rho_a(D_1;q) \cdots \rho_a(D_n;q) \, P_\kappa(x;q,t) =
U_\kappa^{(a)}(x;q,t),
\end{equation}
 which is the $q$-generalization of the exponential
operator formula (\ref{expop}). Replacing $q,t$ by $q^{-1},t^{-1}$
in (\ref{expopu}) and using (\ref{inter}) and (\ref{eE}) gives
\begin{equation}\label{expopv}
{1 \over \rho_a(q\widetilde{D_1};q) \cdots 
\rho_a(q\widetilde{D_n};q)} \, P_\kappa(x;q,t) =
V_\kappa^{(a)}(x;q,t).
\end{equation}
where $\widetilde{D_i}$ denotes the operators  (\ref{qdunkl-def})
with $q,t$ replaced by $q^{-1},t^{-1}$.
We remark that since the $D_i$ and  $\widetilde{D_i}$
are degree lowering operators, only a
finite number of terms in the power series expansion of the $q$-exponential
 operators give a non-zero contribution to (\ref{expopu}) and
(\ref{expopv}).

\subsection{First form of eigenoperator}
With the polynomials $\{U_\kappa^{(a)}(x;q,t)\}$ defined by the 
generating function (\ref{gfmu}), we want to derive the
eigenoperator analogous to (\ref{heo}) which has $\{U_\kappa^{(a)}
(x;q,t)\}$ as eigenfunctions. To present our first representation
of this eigenoperator, let $\tau_i$ denote the $q$-shift operator
of the $i$th variable,
\begin{equation}
\tau_i f(x_1,\dots,x_n) = f(x_1,\dots,x_{i-1},qx_i,x_{i+1},\dots,
x_n),
\end{equation}
let $M_1$ denote the first Macdonald operator,
\begin{equation}\label{defM1}
M_1 := \sum_{i=1}^n A_i(t) \tau_i, \quad
A_i(t) := \prod_{p=1 \atop p \ne i}^n {tx_i - x_p \over x_i - x_p},
\end{equation}
put
\begin{equation}\label{defE0}
E_k :=  \sum_{i=1}^n x^k A_i(t) {\partial \over \partial_q x_i}, 
\quad
{\partial \over \partial_q x_i} :={1 - \tau_i \over (1 - q) x_i}
\end{equation}
and write $\widetilde{M_1}$  for the operator
with $q$ and $t$ replaced by $q^{-1}$ and $t^{-1}$. The required
eigenoperator can be written in terms of $\widetilde{M_1}$ and
$E_0$ (c.f.~the form of the eigenoperator for the generalized
Hermite polynomials given in Proposition 3.2 of ref.~\cite{forr96a}).

\begin{prop}\label{p2.1}
 We have
\begin{equation}\label{eiequ}
{\cal H} \, U_\kappa^{(a)}(x;q,t) = \tilde{e}(\kappa)  U_\kappa^{(a)}(x;q,t)
\end{equation}
where
\begin{equation}\label{form1}
{\cal H} = \widetilde{M}_1 - (1+a)[E_0,\widetilde{M}_1]
+a[E_0,[E_0,\widetilde{M}_1]]
\end{equation}
and $\tilde{e}(\kappa) = \sum_{i=1}^n q^{-\kappa_i} t^{-n+i}$.
\end{prop}

\noindent
{\bf Proof.} \quad We first calculate the action of the Macdonald
operator $\widetilde{M_1}^{(x)}$ (the superscript $(x)$ denotes that 
this operator acts on the $x$-variables only) on the l.h.s.~of
the generating function (\ref{gfmu}). We have
\begin{eqnarray}\lefteqn{
\widetilde{M_1}^{(x)} \prod_{i=1}^n \rho_a(x;q) \,
{}_0 {\cal F}_0(x;y;q,t) } \nonumber \\ &&
= \sum_{p=1}^n A_p^{(x)}(t^{-1}) \rho_a(x;q) \,
 (1-q^{-1} x_p)(1 - q^{-1}ax_p) \tau^{-1}_{p,x} \,
{}_0 {\cal F}_0(x;y;q,t) \nonumber \\ &&
= \prod_{i=1}^n \rho_a(x;q) \Big ( \widetilde{M_1}^{(x)} -
q^{-1}(1+a) \sum_{p=1}^n x_p  A_p^{(x)}(t^{-1})  \tau^{-1}_{p,x}
+ aq^{-2}  \sum_{p=1}^n x_p^2 \nonumber \\
&& \hspace*{2cm} \times A_p^{(x)}(t^{-1})  \tau^{-1}_{p,x}
\Big ) {}_0 {\cal F}_0(x;y;q,t)
\label{arr}
\end{eqnarray}
Using the readily verified identities
\begin{equation}\label{comm1}
[\widetilde{M_1}^{(x)},p_1(x)] =  (q^{-1} - 1)
\sum_{p=1}^n x_p  A_p^{(x)}(t^{-1})  \tau^{-1}_{p,x}
\end{equation}
\begin{equation}\label{comm2}
[[\widetilde{M_1}^{(x)},p_1(x)],p_1(x)] = (q^{-1} - 1)^2
\sum_{p=1}^n x_p^2  A_p^{(x)}(t^{-1})  \tau^{-1}_{p,x}
\end{equation}
where $p_1(x) := \sum_{i=1}^n x_i$, the sums over $p$ in (\ref{arr})
can be replaced in favour of the commutators in 
(\ref{comm1}) and (\ref{comm2}). But from the fact that 
$\widetilde{M}_1$ is an
eigenoperator of $\{P_\kappa(x;q,t)\}$ we see from the definition
of ${}_0 {\cal F}_0$ that
\begin{equation}
\widetilde{M_1}^{(x)} {}_0 {\cal F}_0(x;y;q,t) = 
\widetilde{M_1}^{(y)} {}_0 {\cal F}_0(x;y;q,t).
\end{equation}
Furthermore, from the work of Lassalle \cite{lass97a} we know that
\begin{equation} \label{lassprop}
(1-q) E_0^{(y)} \, {}_0 {\cal F}_0(x;y;q,t) =  p_1(x) 
\, {}_0 {\cal F}_0(x;y;q,t)
\end{equation}
Hence (\ref{arr}) can be rewritten to read
\begin{eqnarray}
\widetilde{M_1}^{(x)} \prod_{i=1}^n \rho_a(x_i;q) \,
{}_0 {\cal F}_0(x;y;q,t)
& = & 
\Big ( \widetilde{M}_1^{(y)} - (1+a)[E_0^{(y)},\widetilde{M}_1^{(y)}]
+a[E_0^{(y)},[E_0^{(y)},\widetilde{M}_1^{(y)}]] \Big ) \nonumber \\
& & \times \prod_{i=1}^n \rho_a(x_i;q) \,
{}_0 {\cal F}_0(x;y;q,t)
\label{arr1}
\end{eqnarray}
Substituting the r.h.s.~of the generating function (\ref{gfmu})
for $\prod_{i=1}^n \rho_a(x_i;q) \,
{}_0 {\cal F}_0(x;y;q,t)$ allows the action of the 
operator to be computed on the l.h.s.~and implies the stated eigenvalue
equation. \hfill$\Box$

\subsection{Second form of eigenoperator}

It is instructive to compare the eigenvalue  equation (\ref{eiequ}) 
with eigenoperator (\ref{form1}) to the
one satisfied by the one variable Al-Salam\&Carlitz polynomials $U^{(a)}_n$.
Let $\tau$ denote the $q$-shift operator in one variable, and 
$D:=x^{-1}(1-\tau)$. Then
\begin{equation} \label{1v-de}
\left(1 - (1+a)D + aD^2 \right)\tau^{-1}\;
U_n^{(a)} = q^{-n}\,U_n^{(a)}
\end{equation}
with a similar equation for $V_n(x;q)$, with $q\rightarrow q^{-1}$
(eq.~(\ref{1v-de}) is equivalent to the hypergeometric-type 
difference equation given on pg.~125 of Askey and Suslov
\cite{askey93a}). This is of a different form to that implied by
(\ref{form1}) with $n=1$.
A multivariable analogue of (\ref{1v-de}) exists, but first some preliminary
results are needed.

\begin{lemma} \label{paloma}

When acting on the space of symmetric functions 
$$
E_0 = \frac{1}{1-q}\;\sum_{i=1}^n D_i \hspace{3cm}
\widetilde{M}_1 = t^{1-n}\sum_{i=1}^n Y_i^{-1}
$$
\end{lemma}
{\bf Proof.}\quad The first identity was given in \cite[Lemma 5.3]{forr97b}. The
proof of the second identity follows a similar line of reasoning, which we
give for completeness. Let
$$
A_{i,m}(t^{-1}) := \prod_{\stackrel{j=1}{j\neq i}}^{m} 
\frac{t^{-1}x_i-x_j}{x_i-x_j}  \hspace{2cm}
\widetilde{M}_1^{(m)} := \sum_{i=1}^m A_{i,m}(t^{-1})\;\tau_i^{-1} 
$$
We shall in fact prove the following stronger result by induction:
\begin{equation} \label{induct.1}
\widetilde{M}_1^{(m)} = t^{1-m}\;\sum_{i=1}^m Y_i^{-1}, \qquad 
m=1,\cdots, n.
\end{equation}
{}From the explicit representation (\ref{yops}) for $Y_i$ (and hence
$Y_i^{-1}$), we have that acting on symmetric functions, $Y_1^{-1}=\tau_1^{-1}$,
so (\ref{induct.1}) is true in the case $m=1$. By induction, the identity 
(\ref{induct.1}) is equivalent to the identity
$$
Y_m^{-1} = t^{m-1} \widetilde{M}_1^{(m)} - t^{m-2} \widetilde{M}_1^{(m-1)}
$$
The action of the operators $Y_m^{-1}$ on symmetric functions can be generated
recursively from $Y_m^{-1} = T_{m-1}\;Y_{m-1}^{-1}$ (which follows from
(\ref{yt1})). It thus suffices to show that the operators
\begin{eqnarray*}
R_m &:=& t^{m-1} \widetilde{M}_1^{(m)} - t^{m-2} \widetilde{M}_1^{(m-1)}\\
&=& t^{m-1} A_{m.m}(t^{-1}) \tau_m^{-1} + t^{m-2}(1-t)\sum_{i=1}^{m-1}
\frac{x_m}{x_i-x_m}\:A_{i,m-1}(t^{-1})\:\tau_i^{-1}
\end{eqnarray*}
obey the same relation: $R_m = T_{m-1}\,R_{m-1}$. This can be verified 
directly. \hfill $\Box$

\bigskip
Using the above lemma, we can now prove the following result

\begin{prop}
Acting on the space of symmetric functions, we have 
\begin{eqnarray}
\left[ E_0, \widetilde{M}_1 \right] &=& \sum_{i=1}^n t^{1-i}
D_i\,Y_i^{-1} \label{cool.1}\\
\left[ E_0, \left[ E_0, \widetilde{M}_1 \right] \right] &=& \sum_{i=1}^n 
t^{1-i} D_i^2\,Y_i^{-1} + (1-t^{-1})\sum_{1\leq i<j\leq n} t^{1-i}
D_j\,D_i\,Y_i^{-1} \label{cool.2}
\end{eqnarray}

\end{prop}
{\bf Proof.}\quad We note from \cite[Lemma 3.3]{forr97b} that
\begin{eqnarray}
\left[Y_j^{-1}, D_i \right] &=& \left\{\begin{array}{cc}
t^{i-j}(1-t) T_{ij}\,D_j\, Y_j^{-1} \qquad i < j \\[2mm]
t^{j-i}(1-t) Y_j^{-1}\,Y_i\,T_{ji}\,D_i\,Y_j^{-1} \qquad i > j  
\end{array} \right.  \label{piel.1}\\
\left[Y_i^{-1}, D_i \right] &=& (q-1)D_i\,Y_i^{-1} + 
(t-1)\sum_{p=i+1}^n t^{-p+i} Y_i^{-1}\,Y_p\,T_{ip}\,D_p\,Y_i^{-1}\nonumber\\
&& +\;q(t-1)\;\sum_{p=1}^{i-1}t^{p-i}\,T_{ip}\, D_i\,Y_i^{-1} \label{piel.2}
\end{eqnarray}
If we now consider the commutator $\left[ E_0, \widetilde{M}_1 \right]$, use
Lemma \ref{paloma} and then (\ref{piel.1}), (\ref{piel.2}), we find that
$$
\left[ E_0, \widetilde{M}_1 \right] = \sum_{i=1}^n D_i\,Y_i^{-1} +(t-1)
\sum_{1\leq i<j\leq n} t^{i-j} T_{ij} D_j Y_j^{-1}
$$
But for $i<j$, 
$T_{ij} D_j Y_j^{-1} = D_i Y_i^{-1}
T_{ij}$,  which puts the operators $T_{ij}$ on the right in the above
expression and allows us to use the fact that when acting on symmetric functions, 
$T_{ij} = t^{2(j-i)-1}$. Simplification of the resulting expression yields
the right hand side of (\ref{cool.1}). The proof of 
(\ref{cool.2}) proceeds
likewise. \hfill $\Box$

\bigskip
The preceding result allows us to give an alternative form of the 
eigenoperator ${\cal H}$ of the Al-Salam\&Carlitz polynomials given in 
(\ref{form1}), namely
\begin{equation} \label{form.2}
{\cal H} = t^{1-n}\sum_{i=1}^n Y_i^{-1} - (1+a) \sum_{i=1}^n t^{1-i}
D_i\,Y_i^{-1} + a \sum_{i=1}^n t^{1-i} D_i^2\,Y_i^{-1} +
a(1-t^{-1})\sum_{1\leq i<j\leq n} t^{1-i} D_j\,D_i\,Y_i^{-1}
\end{equation}
A comparison with (\ref{1v-de}) show that this correctly reproduces the
one-variable result. We also point out that with this choice of 
${\cal H}$, $a=-1$, $q=t^\alpha$, and $x$ replaced by
$\sqrt{2(1-q)} x$, in the limit $q \to 1$ the eigenvalue equation
(\ref{eiequ}) gives the eigenvalue equation for the generalized Hermite
polynomials,
$$
\tilde{H}^{(H)} H_\kappa(x;\alpha) = -2 |\kappa| H_\kappa(x;\alpha),
$$
where $\tilde{H}^{(H)}$ is given by (\ref{heo}).

\subsection{Hermiticity}

In this section it will be shown that for $t=q^{k}$, $k \in\zz^+$,
the polynomials are orthogonal with respect to the inner product
\begin{equation}\label{inneru}
\langle f | g \rangle^{(U)} :=
\int_{[a,1]^n} f(x)g(x) \, d_q\mu^{(U)}(x), \quad
d_q\mu^{(U)}(x) :=   \Delta_q^{(k)}(x)
\prod_{l=1}^n w_U(x_l;q) d_qx_l 
\end{equation}
where
\begin{eqnarray}
\Delta_q^{(k)}(x_1,\dots,x_n) :=
\prod_{p=-(k-1)}^k 
\prod_{1 \le i < j \le n}(x_i - q^px_j)
\end{eqnarray}
and $f$ and $g$ are assumed symmetric.
This will be done by verifying that the eigenoperator ${\cal H}$
in (\ref{eiequ}) is Hermitian with respect to (\ref{inneru}).
Actually we have not been able to deduce the Hermiticity of the 
Al-Salam\&Carlitz eigenoperator ${\cal H}$ directly from either
(\ref{form1}) or (\ref{form.2}). A third
form is required, more amenable to contemplation of this feature. 

First, we write 
$$
D_i = x_i^{-1}\left( 1-t^{2n-i-1}I_{i,n-1}^{-1}\,Y_i \right)
\qquad I_{i,n-1}^{-1}:= T_i^{-1}\cdots T_{n-1}^{-1}\,T_{n-1}^{-1}\cdots
T_i^{-1}
$$
Thus 
\begin{eqnarray}
\sum_{i=1}^n t^{1-i} D_i\,Y_i^{-1} &=& \sum_{i=1}^n t^{1-i} \left(
x_i^{-1}\,Y_i^{-1} - t^{2n-i-1} x_i^{-1}\, I_{i,n-1}^{-1} \right)\nonumber\\
&=& \sum_{i=1}^n t^{-n+i} T_{i-1}^{-1}\cdots T_1^{-1}\,x_1^{-1}\,
\om^{-1}\, T_{n-1}\cdots T_i \label{crest.1}
\end{eqnarray}
and
$D_jD_iY_i^{-1} = A^{(i,j)}_1 + A^{(i,j)}_2 + A^{(i,j)}_3 +
A^{(i,j)}_4$, where
\begin{eqnarray*}
A^{(i,j)}_1 = x_j^{-1} x_i^{-1} Y_i^{-1} &&
A^{(i,j)}_2 = -t^{2n-j-1} x_j^{-1} I_{j,n-1}^{-1} Y_j  x_i^{-1} Y_i^{-1} \\
A^{(i,j)}_3 = -t^{2n-i-1} x_j^{-1} x_i^{-1} I_{i,n-1}^{-1} &&
A^{(i,j)}_4 = t^{4n-i-j-2} x_j^{-1} I_{j,n-1}^{-1} Y_j  x_i^{-1} 
I_{i,n-1}^{-1} 
\end{eqnarray*}
Our aim is to write each of the terms $A^{(i,j)}_p$ so that all occurrences
of $T_i^{\pm 1}$ are on the left of $\omega^{\pm 1}$, while shifting as many
$T_i^{\pm 1}$ as possible to the right 
(whence replacing them by $t^{\pm 1}$ when acting on $f$ or $g$ in
(\ref{inneru})).

Starting with $A_3^{(i,j)}$, we use the fact that acting on symmetric
functions $T_i^{-1} = t^{-1}$, whence
\begin{equation}
A_3^{(i,j)} = -t^{i-1} x_i^{-1} x_j^{-1}, \quad 1\leq i\leq j \leq n,
\label{crest.2}
\end{equation}
which is certainly Hermitian for all $i\leq j$.

Turning to $A_2^{(i,j)}$, simplification with (\ref{id.2}) gives
$$
A_2^{(i,j)} = -q T_j\cdots T_{n-1} \omega T_1\cdots T_{j-1} x_j^{-1}
x_i^{-1} T_{i-1}\cdots T_1\omega^{-1}
$$
To proceed further, we must consider the cases $i=j$, $i<j$ separately.
with the results that
$$
A_2^{(i,j)} = \left\{\begin{array}{ll}
-qt^{i-1} T_j\cdots T_{n-1} \omega T_1\cdots T_i T_{i-1}^{-1}
\cdots T_1^{-1} x_1^{-1} x_i^{-1} \omega^{-1} & i<j \\
-qt^{i-1} T_j\cdots T_{n-1} \omega \left( (t^{-1}-1){\displaystyle 
\sum_{p=2}^i}
T_1\cdots T_{p-1} T_{p-2}^{-1}\cdots T_1^{-1} x_p^{-1} + x_1^{-1}
\right) x_1^{-1} \omega^{-1} & i=j
\end{array} \right.
$$
Combining these together we obtain after further manipulation
\begin{eqnarray}
\sum_{i=1}^nt^{1-i}A_2^{(i,i)} + (1-t^{-1})\sum_{1\leq i<j\leq n}
t^{1-i}A_2^{(i,j)} = q\sum_{i=1}^n T_i\cdots T_{n-1} 
\om x_1^{-2}\om^{-1} \nonumber\\
= q^{-1} \sum_{i=1}^n x_i^{-1}\left( (1-t)\sum_{p=i+1}^n x_p^{-1}
+ x_i^{-1} \right) \label{crest.3}
\end{eqnarray}
which is also Hermitian.

For $A_1^{(i,j)}$ and $A_4^{(i,j)}$, similar considerations hold with 
the result that
\begin{eqnarray}
\sum_{i=1}^nt^{1-i}A_1^{(i,i)} + (1-t^{-1})\sum_{1\leq i<j\leq n}
t^{1-i}A_1^{(i,j)} &=& \sum_{i=1}^n T_{i-1}^{-1}\cdots T_1^{-1} 
x_1^{-2} \om^{-1} \label{crest.4} \\
\sum_{i=1}^nt^{1-i}A_4^{(i,i)} + (1-t^{-1})\sum_{1\leq i<j\leq n}
t^{1-i}A_4^{(i,j)} &=& q\sum_{i=1}^n T_i\cdots T_{n-1} \om x_1^{-2}
\label{crest.5}
\end{eqnarray}
Combining (\ref{crest.1})--(\ref{crest.5}) together, we finally
obtain
\begin{eqnarray}
{\cal H} = f(x) + \sum_{i=1}^n \left( t^{-2n+i+1} T_{n-1}\cdots T_1 \om^{-1}
T_{n-1}\cdots T_i   \right.\hspace{4cm}\nonumber\\
 \left. + t^{-n+i} T_{i-1}^{-1}\cdots T_1^{-1}\, (-(1+a)x_1^{-1}
+ax_1^{-2})\,\om^{-1}\,T_{n-1}\cdots T_i 
 + aq T_i\cdots T_{n-1}\,\om\,x_1^{-2} \right) 
\end{eqnarray}
where $f(x)$ contains no operators $T_i^{\pm 1}$, $\om^{\pm 1}$, and
hence is self-adjoint. To see that this expression is Hermitian with respect
to (\ref{inneru}),
denote the 3 terms in the sum as $\alpha_i$, $\beta_i$ and $\gamma_i$
respectively. Then
\begin{eqnarray*}
(\alpha_i + \beta_i)^{\ast} &=& T_i\cdots T_{n-1}\left\{ t^{-2n+i+1}
(\om^{-1})^{\ast}\,T_1\cdots T_{n-i} \right. \\
&& \left.\;+\; t^{-n+i} (-(1+a)x_1^{-1} +ax_1^{-2}
\,\om^{-1})^{\ast} T_1^{-1}\cdots T_{i-1}^{-1} \right\} \\
&=& t^{-n+1} T_i\cdots T_{n-1} \left( (1-x_1^{-1})(1-ax_1^{-1})\,
\om^{-1}\right)^{\ast} = \gamma_i,
\end{eqnarray*}
where ${}^\ast$ denotes the adjoint with respect to
the inner product (\ref{inneru}), and we have used the fact that
$$
\left((1-x_1^{-1})(1-ax_1^{-1})\,\om^{-1} \right)^{\ast} =
a q t^{n-1} \om\,x_1^{-2}
$$
This latter equation follows from the definition (\ref{omega}) of $\omega$
and the adjoint formulas
$$
(T^{-1}_i)^\ast = T_i^{-1}, \: \:
s_i^\ast = {tx_i - x_{i+1} \over x_i - tx_{i+1}}\,s_i, \: \:\:
\tau_1^\ast = q^{-1}t^{-n+1}
\prod_{l=2}^n {(x_1 - tqx_l) \over (tx_1 - qx_l)}\,
(1-x_1)(1-ax_1) \tau_1^{-1}
$$
obtainable directly from the definitions of $T_i,s_i,\tau_1$
and (\ref{inneru}).

Since $\cal H$ is Hermitian with respect to the inner product
(\ref{inneru}) and $\{U_\kappa^{(a)}\}$ are eigenfunctions of
$\cal H$ with distinct eigenvalues, it follows immediately
that $\{U_\kappa^{(a)}\}$ are orthogonal with respect to 
(\ref{inneru}). Also, by replacing $q$ by $q^{-1}$ and using the
results (\ref{qiuv}) and (\ref{inter}), we see from this result that
$\{V_\kappa^{(a)}\}$ are orthogonal with respect to the inner
product
\begin{equation}
\label{innerv}
\langle f | g \rangle^{(V)} :=
\int_{[1,\infty]^n} f(x)g(x) \, d_q\mu^{(V)}(x), \quad
d_q\mu^{(V)}(x) :=  \Delta_q^{(k)}(x)
\prod_{l=1}^n w_V(x_l;q) d_qx_l .
\end{equation}

\section{Special properties and normalization}
\setcounter{equation}{0}
\subsection{Special properties}
One can check from the generating function (\ref{gfu}) that
the 1-variable polynomials $U_n$ satisfy the formulas
\begin{eqnarray}
\frac{\partial}{\partial_q x}\;U^{(a)}_n(x;q) 
&=& [n]_q\; U_{n-1}^{(a)}(x;q), \qquad \qquad
[n]_q := {1-q^n \over 1 - q} \label{f.1}\\
U_n^{(a/q)}(x;q) &=& U_n^{(a)}(x;q) -q^{-1}a(1-q^n) U_{n-1}^{(a)}(x;q)\\
x\,U^{(a)}_n(x;q) &=& U_{n+1}^{(a)}(x;q) +(1+a)q^n U^{(a)}_n(x;q)
-aq^{n-1}(1-q^n) U_{n-1}^{(a)}(x;q) \label{3term}\\
U_n^{(a)}(1;q) &=& q^{n(n-1)/2}\,(-a)^n \label{f.4}\\
U_n^{(a)}(a;q) &=& q^{n(n-1)/2}\,(-1)^n \label{f.5}\\
U_n^{(0)}(y;q) &=& y^n ({1 \over y};q)_n \label{f.6}
\end{eqnarray}
All these formulas have multivariable generalizations. To
present these generalizations further notions from the theory
of Macdonald polynomials are required. 

First we must introduce multivariable analogues of the $q$-binomial
coefficients $\left ( { \lambda \atop \nu } \right )_{q,t}$
associated with the symmetric Macdonald polynomials, which have
been the subject of a number of recent studies \cite{kaneko96a,okoun97a,%
lass97a,sahi97a}. Lassalle's \cite{lass97a} starting definition takes the form
\begin{equation}\label{las}
P_{\mu}(x;q,t)\;\prod_{i=1}^n \frac{1}{(x_i;q)_{\infty}} =
\sum_{\lambda}\pbin{\lambda}{\mu}_{\!q,t}t^{b(\lambda)-b(\mu)}
\frac{h'_{\mu}}{h'_{\lambda}}\;P_{\lambda}(x;q,t).
\end{equation}
Of particular importance is 
$\left ( { \lambda \atop \nu } \right )_{q,t}$ when $\lambda$ and $\mu$ 
differ by a single node. 
Thus, for a given partition $\lambda$, let $\lambda_{(p)}$ (resp. $\lambda^{(p)}$)
denote the partition $\lambda$ with a node removed from (resp. added to) 
the $p$th row, provided the resulting diagram still represents a valid 
partition.
The corresponding binomial coefficients appear in the following
Macdonald polynomial identities,
\begin{eqnarray}
E_0\:P_{\la}(x;q,t) &=& \sum_{i} \pbin{\la}{\la_{(i)}}_{\!q,t} 
\frac{P_{\la}(t^{\bar{\de}};q,t)}
{P_{\la_{(i)}}(t^{\bar{\de}};q,t)} \, P_{\la_{(i)}}(x;q,t) \label{spider.1}\\
e_1(x)\,P_{\la}(x;q,t) &=& (1-q)\sum_i t^{i-1}\pbin{\la^{(i)}}{\la}_{\!q,t} 
\frac{h'_{\la}}{h'_{\la^{(i)}}}\, P_{\la^{(i)}}(x;q,t) \label{spider.2}
\end{eqnarray}
where $t^{\bar{\de}}:=(1,t,t^2,\ldots,t^{n-1})$, $e_1(x):=\sum_i x_i$, and
the summations are over those $i$ such that $\la_{(i)}$, (resp. $\la^{(i)}$)
is a valid partition.
Furthermore, these particular binomial coefficients have the
explicit evaluation \cite[Theorem 4]{lass97a}
\begin{equation} \label{elem-bin}
\pbin{\lambda}{\lambda_{(p)}}_{\!q,t} = t^{1-p}\frac{1-q^{\la_p}t^{\ell(\la)-p}}
{1-q} \prod_{i=1}^{p-1}\frac{1-q^{\la_i-\la_p}t^{p+1-i}}
{1-q^{\la_i-\la_p}t^{p-i}} \prod_{i=p+1}^{\ell(\la)}
\frac{1-q^{\la_p-\la_i}t^{i-p-1}}{1-q^{\la_p-\la_i}t^{i-p}}
\end{equation}
where $\ell(\la)$ denotes the {\it length} of $\la$. As these
are the only class of 
binomial coefficients which occur in our generalizations,
 this result means that the
terms in our expressions can be computed explicitly.

Also required is the Pieri formula \cite[eq.(6.24)(iv)]{mac}
\begin{equation}
e_r(x) P_\mu(x;q,t) = \sum_{\stackrel{\la}{\mbox{\tiny $\la/\mu$ a vertical
$r$-strip}}}\psi_{\lambda/\mu}'
P_\lambda(x;q,t)
\end{equation}
where $e_r(x)$ denotes the $r$th elementary symmetric function in
the variables $x_1,\dots,x_n$ and
\begin{equation}
 \psi_{\lambda/\mu}' := \prod_{(i,j) \in C_{\lambda/\mu} - 
R_{\lambda/\mu}} {1-q^{\lambda_i - j} t^{\lambda_j' - i + 1} \over
1-q^{\lambda_i - j+1} t^{\lambda_j' - i}}
{1-q^{\mu_i - j+1} t^{\mu_j' - i} \over
1-q^{\mu_i - j} t^{\mu_j' - i + 1}}
\end{equation}
($C_{\lambda/\mu}$ ($R_{\lambda/\mu}$) denotes the union of columns
(rows) which intersect $\lambda - \mu$).

The following result generalizes the first two properties,

\begin{prop}\label{propf}
We have
\begin{eqnarray}
E_0\:U^{(a)}_{\la} &=& \sum_i \pbin{\la}{\la_{(i)}}_{\!q,t} \frac{
P_{\la}(t^{\bar{\de}};q,t)}{P_{\la_{(i)}}(t^{\bar{\de}};q,t)}\,
U^{(a)}_{\la_{(i)}} \label{ovaltine.1}\\
U_{\la}^{(a/q)} &=& \sum_{r=0}^n \left(\frac{-a}{q}\right)^r 
\sum_{\stackrel{\mu}{\mbox{\tiny $\la/\mu$ a vertical $r$-strip}}}
\psi'_{\la/\mu}\frac{f_{\mu}}{f_{\la}}\;U_{\mu}^{(a)} \label{ovaltine.2}
\end{eqnarray}
where $f_{\la}:=t^{b(\la)}/(h'_{\la}(q,t)P_{\la}(t^{\bar{\de}};q,t))$.
\end{prop}
{\bf Proof.}\quad The proof of (\ref{ovaltine.1}) follows a similar
calculation in the Hermite case \cite[Prop. 3.4]{forr96a}, while that of
(\ref{ovaltine.1}) follows a similar result for Laguerre polynomials
\cite[Prop. 4.8]{forr96a}. \hfill $\Box$

\begin{prop} \label{pieri}
We have
\begin{eqnarray*}
e_1\,U_{\la}^{(a)} &=& (1+a)\,e(\la)\,U_{\la}^{(a)} \;+ \;(1-q)\sum_i 
t^{i-1}\pbin{\la^{(i)}}{\la}_{\!q,t}\frac{h'_{\la}}{h'_{\la^{(i)}}}\,
U^{(a)}_{\la^{(i)}} \\
&&\;-\: a\,(1-q)\sum_i q^{\la_i-1} t^{n-i} \pbin{\la}
{\la_{(i)}}_{\,q,t}
\frac{P_{\la}(t^{\bar{\de}};q,t)}{P_{\la_{(i)}}(t^{\bar{\de}};q,t)}
\,U^{(a)}_{\la_{(i)}}
\end{eqnarray*}
where $e(\lambda) = \sum_{i=1}^n q^{\lambda_i} t^{n-i}$.
\end{prop}
{\bf Proof.}\quad The proof follows a similar calculation done in 
\cite{forr96a} for the multivariable Hermite polynomials, using the 
generating function. The product rule for $q$-derivatives complicates
matters slightly, but we get around that as follows: define an operator
$$
B := \sum_{i=1}^n (1-x_i)(1-ax_i)\,A_i(x;t)\,\frac{\partial}{\partial_q x_i}
$$
which has the important properties
\begin{eqnarray}
B^{(x)} \left( \prod_{i=1}^n \rho_a(x_i) \;{}_0{\cal F}_0(x;y;q,t)\right)
&=& \left(B^{(x)} \prod_{i=1}^n \rho_a(x_i) \right) 
{}_0{\cal F}_0(x;y;q,t) \nonumber\\
&&\;+\;\prod_{i=1}^n \rho_a(x_i) \left( E_0^{(x)}\;{}_0{\cal F}_0(x;y;q,t)
\right) \label{b1} \\
B^{(x)} \prod_{i=1}^n \rho_a(x_i) &=& \left( \frac{at^{n-1}}{1-q} e_1(x)
- \frac{(1+a)[n]_t}{1-q} \right) \prod_{i=1}^n \rho_a(x_i) \label{b2}
\end{eqnarray}
where (\ref{b1}) follows from the product rule for $q$-differentiation
and (\ref{b2}) follows directly from computing $\partial/\partial_q x_i
(\rho_a(x_i))$ and using the identities \cite{kaneko96a}
$$
\sum_{i=1}^n A_i(x;t) = [n]_t, \hspace{2cm}
\sum_{i=1}^n x_i A_i(x;t) = t^{n-1} e_1(x).
$$
Returning to the generating function (\ref{gfmu}) and using the
property (\ref{lassprop}), we have (with $f_\lambda$ as given in
Proposition \ref{propf})
\begin{eqnarray*}
\sum_{\la} f_{\la} e_1(y)U^{(a)}_{\la}(y) P_{\la}(x;q,t) &=&
(1-q) \prod_{i=1}^n \rho_a(x_i) \left\{E_0^{(x)} {}_0{\cal F}_0(x;y;q,t)
\right\} \\
&=& \prod_{i=1}^n \rho_a(x_i) \!\left\{ (1-q)B^{(x)} \!-\! at^{n-1}\,e_1(x)
\!+\!(1+a)[n]_t \right\}\! {}_0{\cal F}_0(x;y;q,t) \\
&=& \sum_{\mu} f_{\mu}U_{\mu}^{(a)}(y) \left\{ (1-q)E_0^{(x)} + (1+a)
\left( [n]_t - (1-q)E_1^{(x)}\right) + \right.\\
&&\left.a(1-q)E_2^{(x)} - at^{n-1}e_1(x) \right\}\,P_{\mu}(x;q,t)
\end{eqnarray*}
where, in the second line we have used (\ref{b1}), (\ref{b2}) and
in the third line, we have used the fact that $B:=E_0-(1+a)E_1+aE_2$.
We can now compute the action of the above operators on $P_{\mu}(x;q,t)$
using (\ref{spider.1}), (\ref{spider.2}) and the relation
$$
E_2^{(x)} = \frac{t^{n-1}}{1-q}\,e_1(x) - \frac{1}{1-q}\left[
E_1^{(x)}, e_1(x) \right]
$$
and then, in the resulting equation, compare coefficients of $P_{\la}(x;q,t)$,
to give us the final result. \hfill $\Box$

\bigskip
For the generalization of (\ref{f.4}) and (\ref{f.5}) we have

\begin{prop}
$$
U_{\la}^{(a)}(t^{\bar{\de}};q,t) = (-a)^{|\la|}\,q^{b(\la')}\,t^{-b(\la)}\,
P_{\la}(t^{\bar{\de}};q,t)
\hspace{.7cm} 
U_{\la}^{(a)}(at^{\bar{\de}};q,t) = (-1)^{|\la|}\,q^{b(\la')}\,t^{-b(\la)}
P_{\la}(t^{\bar{\de}};q,t)
$$
\end{prop}
{\bf Proof.}\quad The functions ${}_0{\cal F}_0(x;y;q,t)$ and 
${}_0\psi_0(x;y;q,t)$ obey generalizations of (\ref{e}) and (\ref{E})
in the case $y_i=ct^{i-1}$, (see \cite{kaneko96a,macunp1})
\begin{eqnarray}
{}_0{\cal F}_0(x;ct^{\bar{\de}};q,t) &=& \prod_{i=1}^n 
\frac{1}{(cx_i;q)_{\infty}} \label{tarapoto.1}\\
{}_0\psi_0(x;ct^{\bar{\de}};q,t) &=& \prod_{i=1}^n (cx_i;q)_{\infty}, 
\label{tarapoto.2}
\end{eqnarray}
where it is assumed $|q| < 1$.
Using (\ref{tarapoto.2}) along with the generating function (\ref{gfmu}), 
we have
$$
\sum_{\la}\frac{t^{b(\la)}}{h'_{\la} P_{\la}(t^{\bar{\de}};q,t)}
U^{(a)}_{\la}(t^{\bar{\de}};q,t)P_{\la}(x) = \prod_{i=1}^n (ax_i;q)_{\infty}
\!= {}_0\psi_0(x;at^{\bar{\de}};q,t) 
= \sum_{\la}\frac{(-a)^{|\la|} q^{b(\la')}}{h'_{\la}}\,P_{\la}(x;q,t)
$$
which yields the first formula on comparison of coefficients of
$P_{\la}(x;q,t)$. The second formula follows similarly. \hfill $\Box$

Finally we present the generalization of (\ref{f.6}).

\begin{prop}
Let $P_\lambda^*(x;q,t)$ denote the shifted Macdonald polynomial
\cite{okoun97a,knop96b,sahi97a,lass97a}, which is the unique
(up to normalization)
polynomial of degree $\le|\lambda|$ symmetric in the variables
$x_it^{1-i}$ which vanishes at $(q_1^{\mu_1},\dots,q_n^{\mu_n})$ for
partitions $\mu \ne \lambda$, $|\mu| \le |\lambda|$. With the
normalization such that the leading term of $P_\lambda^*(xt^{-\bar{\delta}}
;q,t)$
is $P_\lambda(x;q,t)$ we have
$$
U_\lambda^{(0)}(y;q,t) = t^{(n-1)|\la|}
P_\lambda^*(yt^{\bar{\delta}-n+1};q,t).
$$
\end{prop}

\noindent
{\bf Proof.} \quad Replace $x$ by $t^{-(n-1)}x$ in  (\ref{las}), 
multiply both sides  by
$$
(-1)^{|\mu|}t^{(n-1)|\mu|}q^{b(\mu')}{P_\mu(y;q,t) \over h_\mu'(q,t)
P_\mu(t^{\bar{\delta}};q,t)}
$$
and sum over $\mu$ using (\ref{kan}) to obtain
\begin{eqnarray*}\lefteqn{
{}_0\psi_0(x;y;q,t) \prod_{i=1}^n{1 \over (t^{-(n-1)}x_i;q)_\infty}}\\
&&= \sum_{\lambda} {t^{-(n-1)|\lambda|}
t^{b(\lambda)} \over h_\lambda'(q,t)}
P_\lambda(x;q,t) \sum_{\mu} \left ( {\lambda \atop \mu} \right )_{\!q,t}
{q^{b(\mu')}t^{(n-1)|\mu|} \over t^{b(\mu)}} (-1)^{|\mu|}
{P_\mu(y;q,t) \over P_\mu(t^{\bar{\delta}};q,t)}.
\end{eqnarray*}
Comparing this with the generating function (\ref{gfmv}) with $a=0$ gives
\begin{equation}\label{vbin}
(-1)^{|\lambda|} q^{b(\lambda')} {V_\lambda^{(0)}(y;q,t) \over
 P_\lambda(t^{\bar{\delta}};q,t)} =
 t^{b(\lambda)} t^{-(n-1)|\la|}
 \sum_{\mu} \left ( {\lambda \atop \mu} \right )_{\!q,t}
 {q^{b(\mu')} \over t^{b(\mu)}} (-1)^{|\mu|} t^{(n-1)|\mu|}
 {P_\mu(y;q,t) \over P_\mu(t^{\bar{\delta}};q,t)}.
\end{equation}
But the shifted Macdonald polynomials have the property that
\cite{okoun97a,lass97a}
\begin{equation}\label{okeq}
{P_\lambda^*(y;q^{-1},t^{-1}) \over P_\lambda^*(0;q^{-1},t^{-1})}
= \sum_{\mu} \left ( {\lambda \atop \mu} \right )_{\!q,t}
 {q^{b(\mu')} \over t^{b(\mu)}} (-1)^{|\mu|}
  {P_\mu(yt^{\bar{\delta}};q,t) \over P_\mu(t^{\bar{\delta}};q,t)}.
\end{equation}
Furthermore \cite{okoun97a}, $P_\lambda^*(0;q^{-1},t^{-1}) =
(-1)^{|\lambda|}  q^{b(\lambda')}  t^{-b(\lambda)}
P_\lambda(t^{-\bar{\delta}};q,t)$, so substituting (\ref{okeq}) in
(\ref{vbin}) we see that
$$
 V_\lambda^{(0)}(y;q,t) = t^{-(n-1)|\lambda|}
 P_\lambda^*(yt^{-\bar{\delta}+n-1};q^{-1},t^{-1}) .
$$
The stated equation follows from this by replacing $q$ and $t$ by their
reciprocals. \hfill$\Box$

\hspace{.2cm}
$\!$We remark that it is possible to use (\ref{vbin}) to deduce an operator
formula relating $V_\lambda^{(0)}(y;q,t)$ and $\left ( {\lambda \atop
\mu} \right )_{q,t}$, which is the symmetric $q$-analogue of a formula
given in ref.~\cite[proof of Prop.~3.23]{forr96d}. For this we require a
$q$-analogue of the Dunkl pairing \cite{dunkl91a}.

\begin{df}
With $D_i$ given by (\ref{qdunkl-def}), and $p$ and $q$ symmetric
polynomials of $n$ variables, define the pairing
\begin{equation}\label{qdp}
[p,q] := p(D) q \Big |_{x=0},
\end{equation}
where $p(D)$ is the operator obtained from $p(x)$ by replacing
each $x_i$ by $D_i$.
\end{df}

For the Dunkl pairing defined by (\ref{qdp}) with the $D_i$ replaced by
the $A$ type Dunkl operators 
$$
d_i := {\partial \over \partial x_i} - \sum_{j \ne i}
{1 - s_{ij} \over x_i - x_j},
$$
the Jack polynomials form an orthogonal set. For the pairing
(\ref{qdunkl-def}) as written, the Macdonald polynomials form an orthogonal
set.

\begin{prop}\label{pfor}
We have
$$
[P_\kappa(x;q,t),P_\sigma(x;q,t)] = t^{-b(\kappa)} h_\kappa'(q,t)
P_\kappa(t^{\bar{\delta}};q,t) \, \delta_{\kappa,\sigma}.
$$
\end{prop}

\noindent
{\bf Proof.} \quad Take $f=P_\kappa(x;q,t)$ in (\ref{fF}) and set $x=0$.
This gives
$$
P_\kappa(D^{(x)};q,t) {}_0{\cal F}_0(x;y;q,t) \Big |_{x=0} =
P_\kappa(y;q,t).
$$
Equating coefficients of $P_\sigma(y;q,t)$ on both sides gives the stated
result.
\hfill$\Box$

\vspace{.2cm}
Acting on both sides of (\ref{vbin}) with $P_\mu(D^{(y)};q,t)$ and setting
$y=0$, Proposition \ref{pfor} immediately implies the sought operator
formula.

\begin{prop}
We have
$$
(-1)^{|\lambda| - |\mu|} {q^{b(\lambda') - b(\mu')} t^{(n-1)(|\lambda| -
|\mu|)}t^{2b(\mu) - b(\lambda)} \over
P_\lambda(t^{\bar{\delta}};q,t) h_\mu'(q,t)}
P_\mu(D^{(y)};q,t) \, V_\lambda^{(0)}(y;q,t) \Big |_{y=0} =
\left ( { \lambda \atop \mu} \right )_{\! q,t}.
$$
\end{prop}

\subsection{Normalization integral for $\lambda = 0$}
For $f=g=1$, the inner product (\ref{inneru}) reads
\begin{equation}\label{nou}
\int_{[a,1]^n} 
\prod_{l=1}^n {E_q(-qx_l) E_q(-{qx_l \over a}) \over
E_q(-q) E_q(-a) E_q(-{q \over a})} 
\prod_{p=-(k-1)}^k 
\prod_{1 \le i < j \le n}(x_i - q^px_j) \, d_qx_1 \cdots d_qx_n
 := {\cal N}_0^{(U)}(a;q,t)
\end{equation}
To evaluate ${\cal N}_0^{(U)}$ we use a modification of the method
used by Aomoto \cite{aomoto94} to evaluate the $q$-Selberg integral.
There are two main steps. The first is to specify the dependence
on $a$. This we do by relating
${\cal N}_0^{(U)}(aq;q,t)$ to ${\cal N}_0^{(U)}(a;q,t)$. In fact
we will show that
\begin{equation}\label{est}
{\cal N}_0^{(U)}(aq;q,t) = t^{n(n-1)/2} {\cal N}_0^{(U)}(a;q,t),
\end{equation}
which implies
\begin{equation}\label{no}
{\cal N}_0^{(U)}(a;q,t) = a^{kn(n-1)/2} g(q,t).
\end{equation}

The second step is to observe that the unknown function $g(q,t)$
can be written as a limit:
\begin{equation}\label{lim}
 g(q,t) = \lim_{a \to 0} a^{-kn(n-1)/2}{\cal N}_0^{(U)}(a;q,t),
\end{equation}
and to evaluate the limit. Inspection of the summand of the 
multiple sum defining ${\cal N}_0^{(U)}(a;q,t)$ for $a \to 0$
shows that the limiting behaviour is determined by the term
$x=t^{\bar{\delta}}$.
Hence
\begin{equation}\label{asym}
{\cal N}_0^{(U)}(a;q,t) \: \sim \: (1-q)^n t^{\sum_{j=1}^n (j-1)}
\prod_{i=1}^n w_U(x_i;q) \Delta_q^{(k)}(x) \Big |_{x = t^{\bar{\delta}}},
\end{equation}
where the notation used in (\ref{inneru}) has been reintroduced.
 Now, for $a \to 0$ and $t=q^k$,
$$
w_U(t^{j-1}q) \: \sim \: q^{-{1 \over 2} k (j-1)(k(j-1)+1)}
{(-a)^{k(j-1)} \over (q;q)_{k(j-1)}}
$$
while
$$
 \Delta_q^{(k)}(t^{\bar{\delta}}) =
t^{2k \sum_{j=1}^{n-2} j (n-1-j)}{(q;q)_{kn}
\Big ( \prod_{l=1}^{n-1} (q;q)_{k(n-l)} \Big )^2 \over
(q;q)_k^n}.
$$
Making use of the summation $\sum_{l=1}^{n-1} l^2 = {1 \over 3} n^3
-{1 \over 2} n^2 + {1 \over 6} n$, substituting these expressions in
(\ref{asym}), then substituting the resulting expression in
(\ref{lim}) gives the evaluation
\begin{equation}\label{qmehta}
{\cal N}_0^{(U)}(a;q,t) =
(1-q)^n (-a)^{kn(n-1)/2}
t^{k \left ({n \atop 3} \right ) -{k-1 \over 2}
\left ({n \atop 2} \right )} 
\prod_{l=1}^n {(q;q)_{kl} \over (q;q)_k}.
\end{equation}
We remark that with $a=-1$, in the limit
$q \to 1$ (\ref{qmehta}) gives a limiting case of the Selberg
integral known as the Mehta integral,
\begin{equation}
\prod_{l=1}^n \int_{-\infty}^\infty d t_l \,
e^{-{1 \over 2} t_l^2} \prod_{1 \le j < k \le n}
|t_k - t_j|^{2k}  = 
(2 \pi)^{n /2} \prod_{j=0}^{n-1} {\Gamma(1+(j+1)k) \over
 \Gamma (1 + k )}
\end{equation}
(this is the integral (\ref{mi}) with $y_j^2 = {1 \over 2} t_j^2$ and $1/\alpha
= k$).
Furthermore, 
by replacing $q$ by $q^{-1}$ and using the result (\ref{qiuv})
we deduce from (\ref{nou}) that for $t=q^k$
\begin{equation}\label{nov}
\int_{[1,\infty]^n} d_q\mu^{(V)}(x) =: {\cal N}_0^{(V)}(a;q,t)
= (1-q)^n a^{kn(n-1)/2}
t^{-2 k \left ({n \atop 3} \right ) -k
\left ({n \atop 2} \right )} 
\prod_{l=1}^n {(q;q)_{kl} \over (q;q)_k}.
\end{equation}

It remains to establish (\ref{est}). Now, it follows immediately
from (\ref{wu}) and (\ref{nou}) that
\begin{equation}\label{insub}
{\cal N}_0^{(U)}(aq;q,t) = 
\int_{[a,1]^n}\, 
\prod_{j=1}^n(x_j - a) \, d_q\mu^{(U)}(x).
\end{equation}
But in Appendix A we will show that
\begin{equation}\label{iden}
\prod_{j=1}^n(x_j - a) = 
\sum_{r=0}^n t^{\bin{n-r}{2}} \; U^{(a)}_{(1^r)}(x)
\end{equation}
Substituting (\ref{iden}) in (\ref{insub}), noting from
(\ref{up}) that $U_{(0)} = 1$, and using the orthogonality of
$\{U_\kappa^{(a)}\}$ with respect to (\ref{inneru}) we obtain
(\ref{est}).

J.~Stokman has pointed out to us that our evaluation
(\ref{qmehta}) is a special case of a result of Evans
\cite{evans94}, in which the ratio of $q$-exponentials
$E_q$ in (\ref{nou}) is replaced by the weight function
(\ref{wbj}) below, and the integration domain
$[a,1]^n$ is replaced by $[-d,c]^n$.

\subsection{Normalization integral for general $\lambda$}

Suppose we have a set of polynomials $\{U_{\la}(x)\}$, orthogonal with
respect to some inner product $\inn{\cdot}{\cdot}$ in which
$e_1$ is self-adjoint. Moreover, suppose
we have a ``first'' Pieri formula of the form
$$
e_1(x)\,U_{\mu}(x) =  \sum_i a_i(\mu)\,U_{\mu^{(i)}}(x) + c(\mu)\,U_{\mu}(x)
+ \sum_i b_i(\mu)\,U_{\mu_{(i)}}(x)
$$
and that there exists a $p$, such that $\mu^{(p)}=\la$. Then using the
orthogonality of the polynomials $\{U_{\la}(x)\}$, we have
\begin{eqnarray}
\inn{U_{\la}}{U_{\la}} &=& \frac{1}{a_p(\mu)}\inn{e_1\,U_{\mu}}{U_{\la}}
= \frac{1}{a_p(\mu)}\inn{U_{\mu}}{e_1\,U_{\la}} \nonumber\\
&=& \frac{b_p(\la)}{a_p(\mu)}\inn{U_{\mu}}{U_{\mu}}
= \frac{b_p(\la)}{a_p(\la_{(p)})}\inn{U_{\la_{(p)}}}{U_{\la_{(p)}}}
\label{rr}
\end{eqnarray}
Now in general, there is no unique way to go from a partition $\nu$ to
a partition $\la\supset\nu$ by adding one node at a time. Indeed, there
are $f^{\la/\nu}$ ways of doing this, where $f^{\la/\nu}$ is the
number of standard tableaux of shape $\la/\nu$. Thus, the relation
(\ref{rr}) can only by iterated {\it provided}
\begin{equation} \label{condition.1}
\frac{b_p(\la)}{a_p(\la_{(p)})} = \frac{g(\la)}{g(\la_{(p)})}
\end{equation}
for some function $g$, in which case there is no inconsistency
i.e.~all sequences going from $\la$ to $\nu$ by removing one node at a 
time will yield the same result. If this condition is satisfied, then clearly
$$
\frac{\inn{U_{\la}}{U_{\la}}}{\inn{U_{\nu}}{U_{\nu}}} =
\frac{g(\la)}{g(\nu)}
$$
In the case of the Al-Salam\&Carlitz polynomials $U^{(a)}_{\la}(x)$,
it follows from Proposition \ref{pieri} that (\ref{condition.1}) is
indeed satisfied with 
$$
g(\la) = (-at^{n-1})^{|\la|} q^{b(\la')} t^{-2b(\la)} h'_{\la}
P_{\la}(t^{\bar{\de}};q,t).
$$
Hence
\begin{equation}\label{normu}
\inn{U_{\la}}{U_{\la}}^{(U)} =: {\cal N}_\lambda^{(U)}(a;q,t)=
(-at^{n-1})^{|\la|} q^{b(\la')} t^{-2b(\la)} h'_{\la}
P_{\la}(t^{\bar{\de}};q,t) {\cal N}_0^{(U)}(a;q,t).
\end{equation}

By replacing $q$ by $q^{-1}$ and using (\ref{qiuv}), (\ref{db}) and
(\ref{ef}) we obtain the normalization formula for $\{V_\kappa^{(a)}\}$
with respect to the inner product (\ref{innerv}),
\begin{equation}\label{normv}
\inn{V_{\la}}{V_{\la}}^{(V)} =: {\cal N}_\lambda^{(V)}(a;q,t)
= (aq^{-1} t^{-2(n-1)})^{|\la|} q^{-2b(\la')} t^{b(\la)} h'_{\la}
P_{\la}(t^{\bar{\de}};q,t) {\cal N}_0^{(V)}(a;q,t).
\end{equation}

\subsection{Integral representations}
The formulas (\ref{normu}) and (\ref{normv}), together with the
orthogonality property of $\{U_\kappa^{(a)}\}$ and
$\{V_\kappa^{(a)}\}$ allow $U_\kappa^{(a)}$ and $V_\kappa^{(a)}$
to be expressed in terms of certain $q$-integrals in which
${}_0^{}{\cal F}_0^{(\alpha)}(x;y)$ and 
${}_0{\psi}_0(x;y;q,t)$ occur as kernels.
Here the convergence properties of 
${}_0{\cal F}_0(x;y;q,t)$ and
${}_0{\psi}_0(x;y;q,t)$ are relevant. Suppose
$0 < q,t <1$. A formula of Macdonald \cite[($7.13'$)]{mac}
gives that $P_\kappa(y;q,t)$ then 
has positive coefficients 
when expressed as a series in monomial symmetric functions, 
so 
we have that $|P_\kappa(y;q,t)| \le
P_\kappa(|y|;q,t)$$ \le c^{|\kappa|} P_\kappa(t^{\bar{\delta}};q,t)$ for
some $c > {\rm max}(t^{-(n-1)} |y_j|)_{j=1,\dots,n}$.
Thus
\begin{equation}\label{bb1}
|{}_0{\cal F}_0(y;x;q,t)| \le
{}_0{\cal F}_0(ct^{\bar{\delta}};|x|;q,t) =
\prod_{i=1}^n {1 \over (c|x_i|;q)_\infty}.
\end{equation}
Similarly we have
\begin{equation}\label{bb2}
|{}_0\psi_0(y;x;q,t)| \le
{}_0\psi_0(-ct^{\bar{\delta}};|x|;q,t) =
\prod_{i=1}^n  (-c|x_i|;q)_\infty.
\end{equation}

The results to be presented are 
the $q$-analogues of Proposition 3.8 and Corollaries 3.1 and
3.2 of ref.~\cite{forr96a}.

\begin{prop}
For $|y|$ and $|z|$ small enough such that all quantities converge,
\begin{eqnarray}\label{p3.8i}\lefteqn{
{1 \over (1-q)^n} \int_{[a,1]^n}{}_0{\cal F}_0(y;x;q,t)
{}_0{\cal F}_0(z;x;q,t) \, d_q\mu^{(U)}(x) \hspace{1.5cm}} \nonumber \\&&
={\cal N}_0^{(U)}(a;q,t) \prod_{l=1}^n {1 \over \rho_a(y_l;q)
 \rho_a(z_l;q)} \, {}_0\psi_0(y;at^{n-1}z;q,t),
\end{eqnarray}
\begin{eqnarray}\label{p3.8ii}\lefteqn{
{1 \over (1-q)^n} \int_{[1,\infty]^n}{}_0{\psi}_0(y;x;q,t)
{}_0{\psi}_0(z;x;q,t) \, d_q\mu^{(V)}(x) \hspace{1.5cm}} \nonumber \\&&
={\cal N}_0^{(V)}(a;q,t) \prod_{l=1}^n {\rho_a(t^{-(n-1)}y_l;q)
 \rho_a(t^{-(n-1)}z_l;q)} \, {}_0{\cal F}_0(y;aq^{-1}t^{-2(n-1)}z;q,t).
\end{eqnarray}
\end{prop}

\noindent
{\bf Proof.} \quad 
The bounds (\ref{bb1}) and (\ref{bb2}), together with the definitions
(\ref{inneru}) and (\ref{innerv}), show that for $|y|$ and $|z|$
small enough the integrals converge.
To verify (\ref{p3.8i}), multiply both sides by
$ \prod_{l=1}^n\! \rho_a(y_l;q)\rho_a(z_l;q)$ and substitute for
$ \prod_{l=1}^n \rho_a(y_l;q) {}_0{\cal F}_0(y;\!x;q,t)$ and 
$ \prod_{l=1}^n \rho_a(z_l;q) {}_0{\cal F}_0(z;\!x;q,t)$ using the
generating function (\ref{gfmu}). Now integrate term-by-term using the
orthogonality of $\{U_\kappa^{(a)}\}$ with respect to the inner product
(\ref{inneru}) and the normalization integral (\ref{normu}), and
identify the resulting series according to the definition (\ref{kan}).
The equation (\ref{p3.8ii}) follows by replacing $q,t$ by $q^{-1},
t^{-1}$ in (\ref{p3.8i}).
\hfill$\Box$

\begin{cor}
We have
\begin{eqnarray}\label{cor1}\lefteqn{
{1 \over (1-q)^n} \int_{[a,1]^n}{}_0{\cal F}_0(y;x;q,t)
U_\kappa^{(a)}(x;q,t) \, d_q\mu^{(U)}(x) \hspace{1.5cm}} \nonumber \\&& =
{\cal N}_0^{(U)}(a;q,t) (-at^{n-1})^{|\kappa|}q^{b(\kappa')}
t^{-b(\kappa)}
 \prod_{l=1}^n {1 \over \rho_a(y_l;q)}
P_\kappa(y;q,t)
\end{eqnarray}
\begin{eqnarray}\label{cor2}\lefteqn{
{1 \over (1-q)^n} \int_{[1,\infty)^n}{}_0{\psi}_0(y;x;q,t)
V_\kappa^{(a)}(x;q,t) \, d_q\mu^{(V)}(x)\hspace{1.5cm}} \nonumber \\&& =
{\cal N}_0^{(V)}(a;q,t) (-aq^{-1}t^{-2(n-1)})^{|\kappa|}q^{-b(\kappa')}
t^{b(\kappa)}
 \prod_{l=1}^n  \rho_a(t^{-(n-1)}y_l;q)
P_\kappa(y;q,t)
\end{eqnarray}
\begin{eqnarray}\label{cor3}\lefteqn{
{1 \over (1-q)^n} \int_{[a,1]^n}{}_0{\cal F}_0(z;x;q,t)
P_\kappa(x;q,t)
\, d_q\mu^{(U)}(x)\hspace{1.5cm}} \nonumber \\&& =
{\cal N}_0^{(U)}(a;q,t) (-t^{n-1})^{|\kappa|}q^{b(\kappa')}
t^{-b(\kappa)}
 \prod_{l=1}^n {1 \over \rho_a(z_l;q)} V_\kappa^{(a)}(az;q,t) 
\end{eqnarray}
\begin{eqnarray}\label{cor4}\lefteqn{
{1 \over (1-q)^n} \int_{[1,\infty)^n}{}_0{\psi}_0(y;x;q,t)
P_\kappa(x;q,t) \, d_q\mu^{(V)}(x)\hspace{1.5cm}} \nonumber \\&& =
{\cal N}_0^{(V)}(a;q,t) (-aq^{-1}t^{-2(n-1)})^{|\kappa|}q^{-b(\kappa')}
t^{b(\kappa)}
 \prod_{l=1}^n \rho_a(t^{-(n-1)}z_l;q) U_\kappa^{(a)}(az;q,t)
\end{eqnarray}
\end{cor}

\noindent
{\bf Proof.} \quad These formulas follow from (\ref{p3.8i}) and (\ref{p3.8ii})
by using the generating functions (\ref{gfmu}) and
(\ref{gfmv}) in an analogous way to the proof of Corollaries 3.1 and 3.2 of
\cite{forr96a}. We remark that (\ref{cor1}) and (\ref{cor2}), and
(\ref{cor3}) and (\ref{cor4}), are equivalent under the mapping 
$q \mapsto q^{-1}$. \hfill$\Box$

\section{Relationship to the big $q$-Jacobi polynomials}
\setcounter{equation}{0}
It has been pointed out to us by J.~Stokman that the
Al-Salam\&Carlitz polynomials $U_\lambda^{(a)}$ introduced herein
can be viewed as the special case $a=b=0$, $c=1$ of his recently
introduced \cite{stok97a} multivariable big $q$-Jacobi
polynomials $P_\lambda^B(\cdot;a,b,c,d;q,t)$, with the parameter
$d$ equal to $-a$ in  $U_\lambda^{(a)}$. This can be seen from the
fact \cite[eq.~(5.18)]{stok97a} that for $t=q^k$, $k \in \zz^+$, the
polynomials $P_\lambda^B$ are orthogonal with respect to the Jackson
integral inner product (\ref{inneru}) modified so that $[a,1]^n$ is
replaced by $[-d,c]^n$ and $w_U(x;q)$ is replaced by
\begin{equation} \label{wbj}
w_B(x;q) := {(qx/c;q)_\infty (-qx/d;q)_\infty \over
(qax/c;q)_\infty (-qbx/d;q)_\infty}.
\end{equation}
The weight function $w_B$ reduces to $w_U$ for the specified
values of the parameters.

An important implication of this identification concerns the orthogonality
for $t \ne q^k$. In \cite{stok97a} it was proved that the $P_\lambda^B$ are
orthogonal with respect to a sum of particular interated Jackson integrals
\cite[eq.~(5.2)]{stok97a}, which reduces to the modification of the
Jackson integral  (\ref{inneru}) specified above when $t=q^k$. For the
choice of parameters specified above in this
more general inner product, the polynomials $U_\lambda^{(a)}$ will thus
form an orthogonal set for general $t \in (0,1)$.

\vspace{.5cm}
\noindent
{\Large\bf Acknowledgements} 

\vspace{3mm}
\noindent
The authors would like to thank Trevor Welsh for helpful remarks,
and J.~Stokman for pointing out the connection with
the big $q$-Jacobi polynomials and associated theory.
This work was supported by the Australian Research Council.

\setcounter{equation}{0}
\setcounter{section}{0}
\renewcommand{\thesection}{\Alph{section}}
\addtocounter{section}{1}
\vspace{3mm}\noindent
{\Large\bf Appendix \thesection}
\vspace{2mm}

In this appendix, we derive the expansion for the product $\prod_{i=1}^n
(x_i-a)$ in terms of Al-Salam--Carlitz polynomials $U^{(a)}_{\lambda}(x)$.
There seems to be no direct way to do this utilizing the generating function 
(\ref{gfmu}), so we proceed indirectly by using the characterization of
$U^{(a)}_{\lambda}(x)$ (resp. $P_{\lambda}(x;q,t)$) as an eigenfunction
of ${\cal H}$ (resp. $\widetilde{M}_1$) in the particular case $\lambda=(1^p)$.
Since the product 
$$
\prod_{i=1}^n (x_i-a)=(-a)^{n}\prod_{i=1}^n(1-x_i/a) = (-a)^{n}
\sum_{r=0}^n e_r(x) \left(\frac{-1}{a}\right)^r
$$
our aim will be to expand the elementary symmetric function $e_r(x)$ in
terms of the $U^{(a)}_{(1^p)}(x)$.
Specifically, we proceed in three stages: Expand $U^{(a)}_{(1^p)}$ in terms
of $e_i$, $0\leq i \leq p$; compute the action of $\widetilde{M}_1$ on
$U^{(a)}_{(1^p)}$ in terms of $U^{(a)}_{(1^i)}$, $i=p$, $p-1$, $p-2$;
expand $e_p$ in terms of $U^{(a)}_{(1^i)}$, $0\leq i \leq p$.

\bigskip\noindent
{\bf Step 1.}\\[2mm]
We first compute the action of the operator ${\cal H}$ on the elementary
symmetric function $e_p$ using the following results:
\begin{eqnarray*}
\widetilde{M}_1\;e_p &=& \tilde{e}(1^p)\,e_p, \hspace{2cm}
\tilde{e}(\lambda) := \sum_{i=1}^n q^{-\lambda_i} t^{-n+i} \\
E_0\;e_p &=& [n+1-p]_t\;e_{p-1}
\end{eqnarray*}
where $[n]_t:=(1-t^n)/(1-t)$. This latter identity appears in the work of
Kaneko \cite{kaneko96a}.
Using these identities and the form of ${\cal H}$ given in Proposition \ref{p2.1},
it follows that 
$$
{\cal H}\;e_p = A_1^{(p)}\,e_p + A_2^{(p)}\;e_{p-1} + A_3^{(p)}\;e_{p-2}
$$
where
\begin{eqnarray*}
A_1^{(p)} &=& \tilde{e}(1^p) \hspace{2cm} A_2^{(p)} = -(1+a)(q^{-1}-1)
t^{-n+p} [n+1-p]_t \\
A_3^{(p)} &=& a(q^{-1}-1)(t-1) t^{-n+p-1} [n+1-p]_t [n+2-p]_t
\end{eqnarray*}
Suppose we have the expansion
$$
U_{(1^p)} = \sum_{i=0}^p a_i^{(p)}\;e_i \hspace{2cm} a^{(p)}_p = 1
$$
If we now use the fact that ${\cal H} U_{(1^p)} = \tilde{e}(1^p) U_{(1^p)}$
and compare coefficients of $e_i$, we obtain a set of 3-term recurrence
relations for the coefficients $a_i^{(p)}$. Indeed, if 
$$
a^{(p)}_{p-i} = f_i(a)\;\qbin{n-p+i}{i}_t , \hspace{2cm}
\qbin{n}{r}_t = \frac{(t;t)_n}{(t;t)_r (t;t)_{n-r}}
$$
then the polynomials $f_i(a)$ obey the 3-term relation
\begin{equation} \label{polys.1}
f_i = -(1+a)\,f_{i-1} + a(t^{i-1}-1)\,f_{i-2} \hspace{2cm}
f_0(a)=1,\qquad f_1(a)=-(1+a)
\end{equation}

Hence
\begin{equation} \label{expansion.1}
U_{(1^p)} = \sum_{i=0}^p f_i(a)\;\qbin{n-p+i}{i}_t\;e_{p-i}.
\end{equation}
By comparison with the 3-term relation for the polynomials $V_n^{(a)}(x;q)$
(see \cite[eq. (5.4)]{ismail85a}), we see that in fact 
$f_i(a):=q^{i(i-1)/2} V_i^{(a)}(0;q)$. 

\bigskip\noindent
{\bf Step 2.}\\[2mm]
We can now use the expansion (\ref{expansion.1}) to calculate the action
of $M_1$ on $U_{(1^p)}$. Indeed, we expect 
\begin{equation} \label{in_bloom}
M_1\,U_{(1^p)} = \gamma_1^{(p)}\, U_{(1^p)} + \gamma_2^{(p)}\, U_{(1^{p-1})}
+ \gamma_3^{(p)}\, U_{(1^{p-2})}
\end{equation}
If we now insert the expansion (\ref{expansion.1}) into this equation and
compare coefficients of $e_{p-i}$, we find that, due to the 3-term relation 
(\ref{polys.1}), 
\begin{eqnarray*}
\gamma_1^{(p)} &=& e(1^p) \hspace{3cm} e(\lambda):=
\sum_{i=1}^n q^{\lambda_i} t^{n-i} \\
\gamma_2^{(p)} &=& -(1-q)(1+a) t^{n-p} [n+1-p]_t \\
\gamma_3^{(p)} &=& (1-q)(t-1)a t^{n-p} [n+1-p]_t [n+2-p]_t
\end{eqnarray*}

\bigskip\noindent
{\bf Step 3.}\\[2mm]
Finally, we can use (\ref{in_bloom}) to obtain the expansion
\begin{equation} \label{stuff.0}
e_p = \sum_{i=0}^p b^{(p)}_i\;U_{(1^i)}, \hspace{2cm} b_p^{(p)} = 1
\end{equation}
This is done by using the fact that $M_1\,e_p = e(1^p)\,e_p$ and following the
procedure used in step 1, to obtain a a 3-term recurrence relation for the
coefficients $b^{(p)}_i$. This recurrence relation has the solution
\begin{equation} \label{stuff.1}
b_{p-i}^{(p)} = \tilde{f}_i(a) \qbin{n-p+i}{i}_t
\end{equation}
where the polynomials $\tilde{f}_i(a)$ obey the 3-term relation
\begin{equation} \label{stuff.2}
\tilde{f}_i = (1+a)\,t^{i-1}\,\tilde{f}_{i-1} + a t^{i-2}(1-t^{i-1})
\tilde{f}_{i-2}, \hspace{2cm} \tilde{f}_0=1, \qquad \tilde{f}_1=1+a
\end{equation}
Again, comparing with the 3-term relation for the polynomials $U_n^{(a)}(x;q)$
given in (\ref{3term}), we see that $\tilde{f}_i(a) = (-1)^i U^{(a)}_i(0;q)$.

\bigskip
We can now use (\ref{stuff.0}), (\ref{stuff.1}) to expand 
$$
\prod_{i=1}^n(1+x_i z) = \sum_{p=0}^n e_p\,z^p = \sum_{r=0}^n
\left(\sum_{j=r}^n b^{(j)}_r\,z^j \right)\;U_{(1^r)} \hspace{2cm}
z:=-1/a
$$
However, the internal summation can be summed exactly, and we claim that
\begin{equation}
\sum_{j=r}^n b^{(j)}_r\,\left(\frac{-1}{a}\right)^j = (-a)^{-n} 
t^{\bin{n-r}{2}}
\end{equation}
Using the explicit expression for the coefficients $b_i^{(p)}$, this is
equivalent to proving (after changing the variables in the summation)
\begin{equation} \label{agradecido}
\sum_{i=0}^{n-r} \tilde{f}_{n-r-i}(a)\qbin{n-r}{i}_t\;(-a)^i = 
t^{\bin{n-r}{2}}
\end{equation}
Denote the left hand side of (\ref{agradecido}) by $S_{n-r}$. It
suffices to show that $S_m = t^{m-1}\,S_{m-1}$. We first note the following
identities for $t$-binomial coefficients:
\begin{eqnarray}
\qbin{m}{i}_t &=& \qbin{m-1}{i}_t + t^{m-i}\qbin{m-1}{i-1}_t \label{gp.1}\\
t^i\qbin{m-1}{i}_t &=& \qbin{m-1}{i}_t - (1-t^{m-i})\qbin{m-1}{i-1}_t\label{gp.2}
\end{eqnarray}
Using (\ref{gp.1}), we have
\begin{eqnarray*}
S_m &=& \sum_{i=0}^m \tilde{f}_{m-i}\left( \qbin{m-1}{i}_t + t^{m-i}
\qbin{m-1}{i-1}_t \right) (-a)^i \\
&=& \sum_{i=0}^{m-1} \tilde{f}_{m-i} \qbin{m-1}{i}_t (-a)^i +
\sum_{i=0}^{m-1} \tilde{f}_{m-i-1}t^{m-i-1}\qbin{m-1}{i}_t (-a)^{i+1}
\end{eqnarray*}
Now apply the 3-term relation (\ref{stuff.2}) to the first term which
yields
\begin{eqnarray*}
S_m &=& \sum_{i=0}^{m-1} \left( t^{m-i-1}\,\tilde{f}_{m-i-1} +
at^{m-i-2}(1-t^{m-i-1})\tilde{f}_{m-i-2} \right) \qbin{m-1}{i}_t (-a)^i \\
&=& \sum_{i=0}^{m-1} t^{m-i-1}\,\tilde{f}_{m-i-1} \left( \qbin{m-1}{i}
-(1-t^{m-i-1})\qbin{m-1}{i-1} \right)\:(-a)^i
\end{eqnarray*}
where to obtain the last line, we have shifted summation variables in 
the second term of the previous line. Using (\ref{gp.2}) thus gives
$S_m = t^{m-1}S_{m-1}$ as desired.

We thus have the expansion (\ref{iden}).

\addtocounter{section}{1}
\setcounter{equation}{0}
\vspace{3mm}\noindent
{\Large\bf Appendix \thesection}
\vspace{2mm}

Here we will show that in the special case $t=q$, the polynomials
$U_\kappa^{(a)}$ admit a determinant formula in terms of their
one-variable counterparts $U_j^{(a)}$. Now, in the case $t=q$,
$\{U_\kappa^{(a)}\}$ form an orthogonal set with respect to the inner
product
$$
\langle f,g\rangle =
\int_{[a,1]^n} \prod_{l=1}^n w_U(x_l;q) 
\prod_{1 \le j < k \le n} (x_j - x_k) (x_j - qx_k) f(x_1,\dots,x_n)
g(x_1,\dots,x_n) \, d_qx_1 \cdots d_qx_n,
$$
where $f$ and $g$ are assumed symmetric.
According to a lemma of Kadell \cite{kad88a}, the product
$\prod_{1 \le j < k \le n} (x_j - qx_k)$ can be replaced by 
$\prod_{1 \le j < k \le n} (x_j - x_k)$ provided the integral is
multiplied by $[n]_q! / n!$. Thus we can write
\begin{equation}\label{innerqq}
\langle f,g\rangle = {[n]_q! \over n!}
\int_{[a,1]^n} \prod_{l=1}^n w_U(x_l;q) 
\prod_{1 \le j < k \le n} (x_j - x_k)^2 f(x_1,\dots,x_n)
g(x_1,\dots,x_n) \, d_qx_1 \cdots d_qx_n,
\end{equation}
and we can characterize the corresponding Al-Salam\&Carlitz polynomials
as the unique symmetric polynomials with an expansion of the form
(\ref{up}) which are orthogonal with respect to (\ref{innerqq}).

Using this characterization, we can easily verify the determinant
formula.

\begin{prop}
We have
\begin{equation}\label{ratio}
U_\kappa^{(a)}(x;q,q) = {\det [U_{\kappa_i + n - i}^{(a)} 
(x_j;q)]_{i,j=1,\dots,n} \over \prod_{1 \le i < j\le n} (x_j - x_i)}
\end{equation}
\end{prop}

\noindent
{\bf Proof.} \quad First note that the ratio (\ref{ratio}) is indeed a
polynomial. The highest degree leading terms are obtained by replacing
each $U_{\kappa_i + n - i}^{(a)}(x_j;q)$ by $x_j^{\kappa_i + n - i}$. This gives the
Schur polynomial $s_\kappa(x)$, and verifies the requirement
(\ref{up}), since $P_\kappa(x;q,q) = s_\kappa(x)$. Now, substituting
(\ref{ratio}) in (\ref{innerqq}) gives
\begin{eqnarray*}
\langle U_\kappa^{(a)}, U_\mu^{(a)} \rangle & = &
{[n]_q! \over n!}
\int_{[a,1]^n} \prod_{l=1}^n w_U(x_l;q) 
\det[U_{\kappa_i + n - i}^{(a)} 
(x_j;q)]_{i,j=1,\dots,n} \det[U_{\mu_i + n - i}^{(a)} 
(x_j;q)]_{i,j=1,\dots,n} \,\\&& \hspace*{2cm} \times d_qx_1 \cdots d_qx_n \\
& = & [n_q]! \int_{[a,1]^n} \prod_{l=1}^n w_U(x_l;q)
U_{\kappa_l + n - l}^{(a)}(x_l;q)\det[U_{\mu_i + n - i}^{(a)} 
(x_j;q)]_{i,j=1,\dots,n} \, d_qx_1 \cdots d_qx_n.
\end{eqnarray*}
Since $\{U_j^{(a)}(x;q)\}$ is an orthogonal set with respect to
the one-dimensional inner product $\langle f, g \rangle :=
\int_a^1 w_U(x;q) f(x) g(x) \, d_qx$, integrating row-by-row in the
determinant  gives
\begin{eqnarray}\label{sun}
\langle U_\kappa^{(a)}, U_\mu^{(a)} \rangle & = &
[n]_q! \Big ( \prod_{i=1}^n \langle U_{\kappa_i + n - i}^{(a)},
U_{\kappa_i + n - i}^{(a)} \rangle \Big ) \delta_{\kappa,\mu} \nonumber \\
& = & [n]_q! (1-q)^n (-a)^{|\kappa|+n(n-1)/2} q^{\sum_{i=1}^n
(\kappa_i +  n - i)(\kappa_i + n - i -1)/2} \prod_{i=1}^n
(q;q)_{\kappa_i + n - i} \, \delta_{\kappa,\mu}, \nonumber \\
\end{eqnarray}
where the final equality follows by using (\ref{ou}), thus establishing
the orthogonality. \hfill$\Box$

\vspace{.2cm}
We remark that straightforward manipulation shows that (\ref{sun}) agrees
with the normalization formulas (\ref{qmehta}) and (\ref{normu}).

\bibliographystyle{plain}

\end{document}